\documentclass[9pt,amstex]{article}
\usepackage{amssymb}
\usepackage{latexsym}

\newtheorem{thm}{Theorem}[section] 
\newtheorem{dfn}[thm]{Definition}
\newtheorem{rmk}[thm]{Remark}
\newtheorem{rmks}[thm]{Remarks}
\newtheorem{cor}[thm]{Corollary}
\newtheorem{prop}[thm]{Proposition}
\newtheorem{lem}[thm]{Lemma}

\newtheorem{ex}[thm]{Example}
\newtheorem{ntt}[thm]{Notations}

\def\tref#1{Theorem~\ref{#1}}
\def\pref#1{Proposition~\ref{#1}}
\def\lref#1{Lemma~\ref{#1}}
\def\cref#1{Corollary~\ref{#1}}
\def\dref#1{Definition~\ref{#1}}
\def\rref#1{Remark~\ref{#1}}
\def\rsref#1{Remarks~\ref{#1}}
\def\rsref#1{Remarks~\ref{#1}}
\def\nref#1{Notations~\ref{#1}}
\def\fref#1{formula~(\ref{#1})}
\def\sref#1{Section~\ref{#1}}

\newcommand{\Pf}{{\em Proof}. }
\newcommand{\EPf}
{%
\mbox{}%
\nolinebreak%
\hfill%
\rule{2mm}{2mm}%
\medbreak%
\par%
}
\newcommand{\C}{\mathbb C} 
\newcommand{\sm}{\star_{\nu}^{M}}
\newcommand{\R}{\mathbb R}

\newcommand{\N}{\mathbb N}
\newcommand{\g}{{\cal G}{}} 
 
\newcommand{\EK}{{\cal L}{}}
\newcommand{\KK}{{\cal K}{}} 

\newcommand{\p}{{\cal P}{}} 
 \newcommand{\e}{{\cal E}{}} 
\newcommand{\pb}{{\overline{\cal P}}{}} 
\newcommand{\s}{{\cal S}{}} 
\renewcommand{\r}{{\cal R}{}}
\renewcommand{\d}{{\cal D}{}} 
\newcommand{\h}{{\cal H}{}} 
\renewcommand{\a}{{\cal A}{}} 
\renewcommand{\b}{{\cal B}{}} 

\newcommand{\z}{{\cal Z}{}}
\newcommand{\Lz}{{\cal Z}{}}
\renewcommand{\k}{{\cal K}{}}

\newcommand{\trianglehorizontal}
{\mbox{\begin{picture}(40,40)\put(0,0){\line(1,0){20}}
\put(-1,-1.5){$\bullet$}
\put(37.5,-1.7){$\bullet$}\put(17.3,27.3){$\bullet$}
\put(40,0){\vector(-1,0){20}}\put(0,0){\vector(2,3){10}}
\put(10,15){\line(2,3){10}}\put(20,30){\vector(2,-3){10}}
\put(30,15){\line(2,-3){10}}\end{picture}}}
\newcommand{\Triangle}[3]{\mbox{\begin{picture}(60,100)
\put(10,20){$#1$}\put(60,20){$#3$}\put(35,60){$#2$}
\put(15,25){\begin{picture}(30,40)\put(0,0){\line(1,0){20}}
\put(-1,-1.5){$\bullet$}\put(37.5,-1.7){$\bullet$}
\put(17.3,27.3){$\bullet$}\put(40,0){\vector(-1,0){20}}
\put(0,0){\vector(2,3){10}}\put(10,15){\line(2,3){10}}
\put(20,30){\vector(2,-3){10}}\put(30,15){\line(2,-3){10}}
\end{picture}}\end{picture}}}

\newcommand{\Trianglefat}[3]{\mbox{\begin{picture}(60,100)
\put(10,20){$#1$}\put(60,20){$#3$}\put(35,60){$#2$}
\put(15,25){\begin{picture}(30,40)\thicklines

\put(0,0){\line(1,0){20}}
\put(-1,-1.5){$\bullet$}\put(37.5,-1.7){$\bullet$}
\put(17.3,27.3){$\bullet$}\put(40,0){\vector(-1,0){20}}
\put(0,0){\vector(2,3){10}}\put(10,15){\line(2,3){10}}
\put(20,30){\vector(2,-3){10}}\put(30,15){\line(2,-3){10}}
\end{picture}}\end{picture}}}

\newcommand{\carrehoriz}[1]{\mbox{\begin{picture}(100,100)
\put(10,10){\line(1,0){70}}\put(10,80){\line(1,0){70}}\put(10,10){\line(1,1){70}}

\put(10,80){\line(1,-1){70}}\put(30,30){\vector(1,1){2}}
\put(60,60){\vector(-1,-1){2}}\put(45,10){\vector(-1,0){2}}
\put(45,80){\vector(1,0){2}}\put(60,30){\vector(1,-1){2}}
\put(30,60){\vector(-1,1){2}}\put(8,8){$\bullet$}\put(8,78){$\bullet$}
\put(77,8){$\bullet$}\put(78,78){$\bullet$}\put(42,42){$\bullet$}\put(0,5){$d$}

\put(85,5){$c$}\put(85,80){$b$}\put(0,80){$a$}\put(50,43){#1}\end{picture}}}

\newcommand{\carrevertic}[1]{\mbox{\begin{picture}(100,100)
\put(10,10){\line(0,1){70}}\put(80,10){\line(0,1){70}}
\put(10,10){\line(1,1){70}}\put(10,80){\line(1,-1){70}}
\put(30,30){\vector(-1,-1){2}}\put(60,60){\vector(1,1){2}}
\put(10,45){\vector(0,1){2}}\put(80,45){\vector(0,-1){2}}
\put(60,30){\vector(-1,1){2}}\put(30,60){\vector(1,-1){2}}
\put(8,8){$\bullet$}\put(8,78){$\bullet$}\put(77,8){$\bullet$}
\put(78,78){$\bullet$}\put(42,42){$\bullet$}\put(0,5){$d$}\put(85,5){$c$}

\put(85,80){$b$}\put(0,80){$a$}\put(50,43){#1}\end{picture}}}
\newcommand{\carreverticfat}[1]{\mbox{\begin{picture}(100,100)
\thicklines
\put(10,10){\line(0,1){70}}\put(80,10){\line(0,1){70}}
\put(10,10){\line(1,1){70}}\put(10,80){\line(1,-1){70}}
\put(30,30){\vector(-1,-1){2}}\put(60,60){\vector(1,1){2}}
\put(10,45){\vector(0,1){2}}\put(80,45){\vector(0,-1){2}}
\put(60,30){\vector(-1,1){2}}\put(30,60){\vector(1,-1){2}}
\put(8,8){$\bullet$}\put(8,78){$\bullet$}\put(77,8){$\bullet$}
\put(78,78){$\bullet$}\put(42,42){$\bullet$}\put(0,5){$d$}\put(85,5){$c$}

\put(85,80){$b$}
\put(0,80){$a$}
\put(50,43){#1}
\end{picture}}}

\begin{document}
\title{Strict Quantization of Solvable Symmetric Spaces}
\author{{\bf Pierre Bieliavsky}\\
D\'epartement de Math\'ematiques\\
Universit\'e Libre de Bruxelles\\
CP. 218, Campus Plaine\\
B-1050, Brussels\\
Belgium\\{email: pbiel@ulb.ac.be}\\}
\date{}
\maketitle
{\center{\sl To the memory of Mosh\'e Flato}}
\begin{abstract}
This work is a contribution to the area of Strict Quantization (in the 
sense of Rieffel) in the presence of curvature and non-Abelian group 
actions. More precisely, we use geometry to obtain explicit oscillatory integral 
formulae for strongly invariant strict deformation quantizations of a class 
of solvable symplectic symmetric spaces. 
Each of these quantizations gives rise to a field of (pre)-$C^{\star}$-algebras 
whose fibers are function algebras which are closed under the 
deformed product. The symmetry group of the symmetric space acts on 
each fiber by $C^{\star}$-algebra automorphisms.
\end{abstract}
\section*{Introduction}\label{INTRO}

Weyl's method for quantizing a free particle in $\R^n$ consists in 
a correspondence between classical observables (i.e.~functions 
on the phase space $\R^{2n}=T^\star(\R^n)$) and linear operators 
acting on the Hilbert space of square integrable functions on the 
configuration space $\R^{n}$. By reading the operator 
composition product at the level of functions 
via Weyl's correspondence, one gets a non-commutative associative 
product on the classical observables called Weyl's product. This new 
product appears as a deformation with one parameter, namely Planck's ``constant" 
$\hbar$, of the usual commutative pointwize multiplication of functions 
in the direction of the classical Poisson bracket. This is the starting 
point of Formal Deformation Quantization (FDQ) theory (or star products), 
as introduced  by Bayen, Flato, Fronsdal, Lichnerowicz and Sternheimer in 
the late seventies \cite{BFFLS}. In FDQ, Quantum Mechanics is formulated 
in a classical framework as the space of formal power series in $\hbar$ 
with classical observables as coefficients. No reference to a particular 
Hilbert space representation is made.  

An important feature of Weyl's product is that it is ``strict",
meaning that in a suitable functional framework, the product 
of two functions is again a function, rather than a formal power 
series in $\hbar$, as in star product theory. For instance,
the Schwartz space on $\R^{2n}$ is closed under Weyl's product.
Originally formulated in the framework of $C^\star$-algebras, the notion 
(theory) of Strict Deformation Quantization (SDQ) was  introduced 
by M. A. Rieffel in 1989 \cite{R2} (see also \cite{N}), and has been extensively 
developed in various directions since then. For example, SDQ provides 
a very elegant and simple description of Quantum Tori by generalizing Weyl's 
quantization of $\R^{2n}$ to Poisson manifolds whose Poisson structure 
comes from an action of $\R^d$ \cite{R1}. In the latter, invariance of 
Weyl's product under translations plays a crucial role.

There is a very simple oscillatory integral formula for Weyl's 
deformed product (see~formula~(\ref{WEYL}) in the next section). 
Such an oscillatory integral formula is not only convenient 
when studying the functional analysis and strictness of the quantization, 
it also indicates a possible approach to the study of the geometry of the 
quantization.

On a symplectic manifold $M$, it is natural to express a 
given (continuous) multiplication, $\star$, on functions via a kernel 
formula of the type~:
$$
u\star v(x)=\int_{M\times M}K(x,\, .\, , \,. \,)u\otimes v ,
$$
where $u$ and $v$ are two functions on	$M$, where $x$ is a point in 
$M$, and where $K$ is the three-point kernel defining the multiplication 
$\star$. The integration is with respect to the product
Liouville measure on $M\times M$.

In the so-called {\sl WKB-quantization} program (see~\cite{W1}) initiated 
independently by Karasev, Weinstein and Zakrzewski, one considers kernels 
of the form~: 
$$
K=a_\hbar\, e^{\frac{i}{\hbar}S}.
$$
Here $S$ is a real-valued smooth function on $M\times M \times M$ 
called {\sl the phase} and  $a_\hbar$, called {\sl the amplitude},
is usually a power series in $\hbar$. It is defined in such a way that the 
product $\star$ would be, at least formally, associative. It should 
also constitute a one-parameter ($\hbar$) deformation of the usual 
pointwize product in the direction of the symplectic Poisson bracket.
Briefly, a WKB-quantization is defined by an oscillatory integral 
product formula.

When one imposes compatibility with some geometric structure given 
on the symplectic manifold, for instance invariance of the 
quantization under the action of the group of affine transformations 
of a given symplectic connection, $S$ and $a_\hbar$ become extremely constrained 
and carry  non-trivial geometric information. In 
the simplest example of Weyl's quantization of $M=\R^{2n}$, the 
amplitude function is identically equal to 1 while the phase 
$S(x,y,z)$ is proportional to the symplectic area of the Euclidean 
triangle ${\Delta}{xyz}$ whose vertices are points $x,y$ and $z$ 
\cite{K2,W1}. 
In particular, in this case the largest group of diffeomorphisms of 
$\R^{2n}$ preserving the quantization is $Sp(n,\R)\times\R^{2n}$, that 
is, the group of symplectic displacements with respect to the 
{\sl flat} connection on $\R^{2n}$.

In this work, we study WKB-quantizations of a class of {\sl curved} affine 
symplectic manifolds. More precisely, a {\bf symplectic symmetric space} 
is a triple $(M,\omega,\nabla)$ where $(M,\omega)$ is a symplectic	manifold 
and where $\nabla$ is a \def\fref#1{formula~\ref{#1}} torsion-free affine 
connection on $M$ such that $\nabla\omega=0$ and such that at every point 
$x\in M$, the local geodesic symmetry, $s_x$, extends globally to $M$ as 
an affine symplectic transformation. In this case the group $G(M)$ 
of transformations of $M$ generated by the symmetries $\{s_x\}$ completely 
determines the connection (see~Section~\ref{SSS}). By a WKB-quantization of 
a symplectic symmetric space, we essentially mean a WKB-quantization of 
$(M,\omega)$ with kernel 
$K=a_\hbar\, e^{\frac{i}{\hbar}S}$ as above such that  

\begin{enumerate}
\item[(i)] the amplitude $a_\hbar$ and the phase $S$ are invariant under 
the group $G(M)$; 
\item[(ii)] for all value of $\hbar\in[0,\infty)$, there exists a function 
algebra $\e_\hbar\subset\mbox{Fun}(M)$ stable under the deformed 
multiplication and containing the space of smooth compactly supported 
functions. 
(We denote by $\mbox{Fun}(M)$ the space of complex 
valued functions on $M$).
\end{enumerate}
For a precise definition see Definition~\ref{WKB}.

One interest of such an invariant quantization is that, when well
behaved, it leads to deformations of quotient spaces 
$\Gamma\backslash M\mbox{ where }\Gamma\subset G(M)$ (see 
Section~\ref{CONCL}). However, this 
last question will not be investigated in the present article.

WKB-quantization of symplectic symmetric spaces has already been 
investigated by A. Weinstein. In \cite{W1}, assuming the existence 
of a WKB-quantization of a symmetric symplectic space which is 
invariant under the group of automorphisms, Weinstein gives a beautiful 
geometric description of the phase $S(x,y,z)$ in terms of the symplectic 
area of a geodesic triangle admitting points $x,y$ and $z$ as midpoints 
of its edges. Even for symmetric spaces, the problems of finding the 
amplitude $a_\hbar$ as well as that of giving a suitable functional framework 
where the WKB-quantization would yield topological function algebras are 
still wide open. Nevertheless, a geometrical study combined with  
techniques 
coming from Star Representation theory \cite{AC2,Fr} (representation 
theory of Lie groups in the framework of deformation quantization) 
seem to suggest a way to attend these questions. At least for solvable 
symmetric spaces, one can give a quite satisfying answer. This is what 
is done in this paper. The symplectic spaces considered here are symplectomorphic 
to $\R^{2n}$ endowed with its standard symplectic structure. From the topological 
point of view the situation is therefore trivial. In contrast, the geometry 
is not. Indeed, each of these spaces is endowed with a curved symplectic affine 
connection whose symplectic affine transformations are not linear. In particular, 
Weyl's product is not preserved under such canonical transformations. 
In our case, the automorphism group $G(M)$ is a solvable Lie group. Our situation therefore 
differs from the K\"ahlerian symmetric case whose study was initiated 
by Berezin \cite{B2} (see Section~\ref{CONCL}).

The present	paper is organized as follows. 
\begin{enumerate}
\item[]{\bf \ref{WEYLR} Weyl's Quantization revisited}\\
We introduce a geometric setting in which 
associativity of Weyl's product is interpreted in elementary 
geometric terms. We propose a convenient 
framework for WKB-quantization (Definition~\ref{WKB}) inspired by
Rieffel's definition of strict deformation quantization. 
\item[] {\bf \ref{SSS} General facts about symplectic symmetric 
spaces}\\
For the convenience of the reader, we recall some relevant definitions 
and results on symplectic symmetric spaces which can be found in 
\cite{Bi3}.
\item[] {\bf \ref{PHSSS} Phase functions on symplectic symmetric spaces}\\
We adapt the geometric setting introduced in 
Section~\ref{WEYLR} to the case of symplectic symmetric spaces. We indicate how Weinstein's phase function 
\cite{W1} appears in this picture.
\item[] {\bf \ref{ESSS} Elementary solvable symplectic symmetric spaces}\\
We define the class of symplectic symmetric spaces that will 
be considered in this work. Those are such that the action of the 
holonomy algebra at a point $o$ has an isotropic range in the 
symplectic tangent space at $o$. The symmetry group of such a symplectic symmetric space 
is a solvable Lie group. Among solvable symplectic symmetric spaces, 
these spaces are structurally of primary importance (see Section \ref{ESSS} and 
\cite{Bi1}). We call them elementary solvable symplectic symmetric spaces. We end this 
section by giving explicit formulae for the invariant phase functions $S$ on 
these spaces as well as an embedding property (Proposition~\ref{SESET}) 
which will be useful when studying the functional analysis attached 
to the integral kernel $K=a_\hbar e^{\frac{i}{\hbar}S}$.
\item[] {\bf \ref{REPTH} Obtaining the oscillating kernel via star 
representation theoretical 
methods}\\
On a symplectic symmetric space as considered in Section \ref{ESSS}, one has a 
canonical Darboux chart. In this Darboux chart, the formal Moyal star 
product	 turns out to be	covariant (in Arnal's sense) under the 
action of
some central extension of the transvection algebra. Using 
the cocycle defining the associated 
star representation, we introduce some kind of 
integral Fourier operator on the quantum algebra which 
intertwines the usual commutative product with an
$\hbar$-dependent {\sl commutative} product for 
which the extended transvection algebra acts by derivations. 
The ``commutative manifold'' underlying the latter product 
carries therefore a strongly invariant deformation
quantization. This space turns out to be equivariantly isomorphic	
to the symmetric space at hand.
\item[] {\bf \ref{WKBQ} WKB-Quantization}\\
We determine explicitly the amplitude functions $a_{\hbar}$	that yield invariant 
WKB-quantizations of our solvable symplectic symmetric spaces via their associated kernels
$K=a_{\hbar}e^{\frac{i}{\hbar}S}$. 
Combined with 
Section~\ref{REPTH}, this yields explicit	oscillatory integral formulae for strongly invariant strict deformation 
quantizations (Theorem \ref{WKBTHM}). We define deformed function algebras
i.e.~we exhibit spaces of functions on the symmetric space which are stable under the 
oscillatory integral deformed product.
\item[] {\bf \ref{TA} Topological algebras}\\
We define symmetry invariant $C^{\star}$-norms on the above mentioned function 
algebras.
\item[] {\bf \ref{CONCL} Remarks for further developments}\\
The present work gives us a hope for attending ``quantum Anosov 
property'' on compact solv-manifolds. Also, we believe it could allow 
to attend ``quantum surfaces'' (\cite{Ra}, \cite{NN}) in an explicit 
and quite elementary way.
\end{enumerate}

\noindent{\bf Acknowledgments}\\
The author would like to express his gratitude to Alan Weinstein	
for many enlightening discussions as well as an invitation in January 
of 2000 at the Mathematics Department of UC Berkeley. Considerable
improvements to this paper were made there. He also would like to thank 
Marc Rieffel for  communicating his work on 
deformation quantization for actions of $\R^d$ and pointing out 
the right framework for ``convergent" deformation quantization.
The author is indebted to M\'elanie Bertelson for having
suggested the important choice of a norm in the proof of Lemma~\ref{GROWTH}.
At last, the author thanks Jeff Kiralis and Stefan Waldmann for reading the 
first draft of the manuscript and suggesting improvements to it. 
This research has been partially supported by the Communaut\'e Fran\c caise de 
Belgique, through an Action de Recherche Concert\'ee de la Direction 
de la Recherche Scientifique.

\section{Weyl's Quantization revisited${}^{\mbox{\footnotesize{1}}}$}
\label{WEYLR}
\footnotetext[1]{In this section 
we will not consider the functional analytical problems involved in existence 
and well-definedness of our oscillatory integrals, those will be investigated 
in the next sections. Only the geometric aspects will be discussed in the 
present section.}

Weyl's quantization consists in a correspondence between
classical observables of a mechanical system, i.e. functions
on a phase space and operators on a suitable Hilbert space.
More precisely, when the phase space is the symplectic
vector space $(\R^{2n},\omega^0)$, Weyl's quantization map
is given by~: 
$$
{\cal S}(\R^{2n}) \stackrel{W_\hbar}{\to} \b(L^2(\R^{n}))
$$
$$
(W_\hbar(u)f)(q) \stackrel{\mbox{}}{=} \int_{\R^{n}\times\R^{n}} 
e^{\frac{i}{\hbar}<q-\eta,\xi>}u\left(\frac{q+\eta}{2},\xi\right) 
\, f(\eta)\, d\xi\, d\eta,
$$
where $\R^{2n}$ is seen as $\R^{n}\times \R^{n}=\{ (q,p)\} $,  $d\xi$  
(respectively $d\eta$) is some suitable normalization of the Haar measure 
on the Abelian Lie group $\R^{n}$, where ${\cal S}(\R^{2n})$ denotes the 
space of Schwartz's functions on $\R^{2n}$ and $\b(L^2(\R^{n}))$ is the 
algebra of bounded operators on $L^2(\R^{n})$. The Weyl product of two 
Schwartz's functions is formally defined as~: 
$$
W_\hbar( u \star^0_\hbar v )= W_\hbar(u) \circ W_\hbar(v).
$$
One derives the following integral formula for the  Weyl product~:
\begin{equation}\label{WEYL}
(u \star^0_\hbar v)(x)=\frac{1}{\hbar^{2n}}\int _{\R^{2n}\times\R^{2n}} 
u(y)\, v(z)\, e^{-\frac{2i}{\hbar}S^{0}(x,y,z)} dy\, dz,
\end{equation}
where $$S^{0}(x,y,z)=\omega^0(x,y)+\omega^0(y,z)+\omega^0(z,x),$$ and 
where  $dy$ (respectively $dz$) is some suitable normalization of the Liouville 
measure.

Endowed with this product, the Schwartz space $({\cal 
S}(\R^{2n}),\star^0_\hbar)$ becomes an  associative topological algebra 
\cite{H}. Its algebra structure extends to $L^2(\R^{2n})$.

Interpreting formula (\ref{WEYL}) as an oscillatory integral with parameter 
$-\frac{2}{\hbar}$, one can use  a stationary phase method to 
obtain the following asymptotic expansion~:

\begin{equation}\label{MOYAL}
u \star^0_\hbar v \sim
 \, uv+\nu\{u,v\}+\sum_{k=2}^{\infty}\frac{\nu^{k}}{k!}\sum_{\begin{array}{c} 
i_{1}\ldots i_{k}\\ j_{1}\ldots j_{k}
\end{array}} \Omega^{i_{1}j_{1}}\ldots \Omega^{i_{k}j_{k}} 
\partial_{i_{1}\ldots i_{k}}u.\partial_{j_{1}\ldots j_{k}}v
\end{equation}
with $\nu=\frac{\hbar}{2i}$ and where 
$\Omega^{ij}\partial_{i}\wedge\partial_{j}=\{\, , \,\}$ is the Poisson  
tensor associated to $\omega^{0}$. The RHS of (\ref{MOYAL}) extends 
by $\C[[\hbar]]$-bilinearity as an associative product to the space
of formal power series $C^\infty(\R^{2n})[[\hbar]]$. 
This formal product is called the Moyal star 
product and will be denoted by $u\star^M_{\nu}v$.

Dirac's condition in Quantum Mechanics in this context reads as
$$
\frac{1}{2\nu}\left(u\star^M_{\nu}v-v\star^M_{\nu}u\right)\,\mbox{mod}(\nu)
=\{u,v\}.
$$ 
The oscillatory integral formula~(\ref{WEYL}) indicates that the triple 
$(\R^{2n},dx, S^{0})$ contains all the structure needed to produce an 
associative non commutative deformation of the usual pointwize product 
on functions on $\R^{2n}$. Of course, compatibility  between  $S^{0}$ and 
the (Abelian) group structure of $\R^{2n}$ is certainly crucial. So, a 
natural question is~: {\sl given an orientable manifold $M$ endowed with 
a volume form $\mu$, what are the  conditions on a three-point function 
$S\in C^{\infty}(M\times M\times M,\R)$ which would guarantee associativity 
of the product 
$$
u\star v(x)=\int_{M\times M}e^{iS(x\, , . \, , \, . \,)}u\otimes v \,
\mu\otimes \mu\; ?
$$}

When writing a (continuous) multiplication on functions via a kernel formula 
of the type~:  
$$
u\star v(x)=\int_{M\times M}K(x,\, .\, , \,. \,)u\otimes v\quad , 
$$
a computation shows that associativity for the multiplication $\star$ is 
(at least formally) equivalent to the following condition~: 

\begin{equation}\label{ASSK}
\int_{M}K(a,b,t)K(t,c,d)\mu(t)=\int_{M}K(a,\tau,d)K(\tau,b,c)\mu(\tau),
\end{equation}
for every quadruple of points $a,b,c,d$ in $M$. Equality~(\ref{ASSK}) obviously 
holds if one can pass from one integrand to the other using a change a 
variables $\tau=\varphi(t)$. This motivates 

\begin{dfn}\label{coucou}
Let $(M,\mu)$ be an orientable manifold endowed with a volume form $\mu$. 
A three-point kernel $K\in C^{\infty}(M\times M\times M)$ is {\bf geometrically 
associative} if for every quadruple of points $a,b,c,d$ in $M$ there exists 
a volume preserving diffeomorphism
$$
\varphi:(M,\mu)\to(M,\mu),
$$
such that for all $t$ in $M$~:
$$
K(a,b,t)K(t,c,d)=K(a,\varphi(t),d)K(\varphi(t),b,c).
$$
\end{dfn}

We will prove in Proposition~\ref{GEOM} that the following structure leads to a 
geometrically associative kernel.
\begin{dfn}\label{WT}
A {\bf Weyl triple} is a triple $(M,\mu,S)$ where $(M,\mu)$ is an oriented 
manifold and $S:M\times M\times M \to {\R}$ is a smooth three-point function 
on $M$ such that~:
\begin{enumerate}
\item[(i)] $S(x,y,z)=S(z,x,y)=-S(y,x,z)$;
\item[(ii)] for all $m \in M$, one has~:
$$
S(x,y,z)=S(x,y,m)+S(y,z,m)+S(z,x,m)\qquad \forall x,y,z\in M;
$$
\item[(iii)] for all $x \in M$, there exists a $\mu$-preserving 
diffeomorphism $s_x:M\to M$ such that~:
$$
S(x,y,z)=-S(x,s_x(y),z)\qquad \forall y,z\in M.
$$
\end{enumerate}
\end{dfn}
Property (i) in Definition~\ref{WT} naturally leads us to adopt the following 
``oriented graph" type notation for~$S$~:

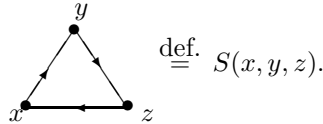
$$\begin{picture}(200,50)
\put(0,0){$x$}
\put(50,0){$z$}
\put(25,40){$y$}
\put(5,5){\trianglehorizontal}
\put(57,19){$\stackrel{\mbox{def.}}{=}\,S(x,y,z).$}
\end{picture}$$
A change of orientation in such an ``oriented triangle'' leads to a change 
of sign of its value. However, the value represented by such a ``triangle''  
does not depend on the way it ``stands'', only the data of the vertices 
and the orientation of the edges matters.

Now, consider a Weyl triple $(M,\mu,S)$, and let $A$ be some (topological) 
associative algebra. And, for compactly supported functions $u$ and $v\in 
C^{\infty}_{c}(M,A)$, consider the following ``product''~:  
$$
u\star v(x)=\int_{M\times M} u(y) v(z) e^{iS(x,y,z)}\mu(y) \,\mu(z).
$$

With the above notation for $S$, associativity for $\star$ now formally reads as follows~:

$$\begin{array}{l}
\begin{array}{ccccc}
\begin{picture}(60,50)
\put(0,25){
$u\star(v\star w)(a)=$}
\end{picture}
&
\begin{picture}(60,80)
\put(0,25){\mbox{$\displaystyle \int \exp i ($}}

\put(30,10){
\begin{picture}(50,35)
\put(0,0){$d$}
\put(30,15){$t$}
\put(0,30){$a$}
\put(5,5){$\bullet$}
\put(5,25){$\bullet$}
\put(25,15){$\bullet$}
\put(7,17){\vector(0,1){2}}
\put(17,22){\vector(2,-1){2}}
\put(17,12){\vector(-2,-1){2}}
\put(7,7){\line(0,1){20}}
\put(7,7){\line(2,1){20}}
\put(27,17){\line(-2,1){20}}
\end{picture}
}
\end{picture}
&

\begin{picture}(60,10)
\put(0,25){\mbox{$) \displaystyle \int \exp i ($}}

\put(40,10){
\begin{picture}(50,35)
\put(25,0){$c$}
\put(-2,15){$t$}
\put(25,30){$b$}
\put(20,5){$\bullet$}
\put(20,25){$\bullet$}
\put(2,15){$\bullet$}
\put(23,17){\vector(0,-1){2}}
\put(13,22){\vector(2,1){2}}
\put(13,12){\vector(-2,1){2}}
\put(23,7){\line(0,1){20}}
\put(23,7){\line(-2,1){20}}
\put(3,17){\line(2,1){20}}
\end{picture}
}
\end{picture}
&
\begin{picture}(70,50)
\put(0,25){
$)\;u(b)\,v(c)\,w(d)\;$}
\end{picture}
&
\end{array} \\

\begin{array}{cccc}
\begin{picture}(60,10)
\put(-7,25){\mbox{$=\displaystyle \int \exp i ($}}

\put(43,15){

\begin{picture}(50,35)
\put(22,25){$b$}
\put(10,-2){$\tau$}
\put(-8,25){$a$}
\put(16,20){$\bullet$}
\put(-4,20){$\bullet$}
\put(6,2){$\bullet$}
\put(8,23){\vector(1,0){2}}
\put(3,13){\vector(-1,2){2}}
\put(13,13){\vector(-1,-2){2}}
\put(18,23){\line(-1,0){20}}
\put(18,23){\line(-1,-2){10}}
\put(8,3){\line(-1,2){10}}
\end{picture}
}
\end{picture}
&

\begin{picture}(55,80)
\put(0,25){\mbox{$) \displaystyle \int \exp i ($}}

\put(43,12){

\begin{picture}(60,35)
\put(21,1){$c$}
\put(6,30){$\tau$}
\put(-10,1){$d$}
\put(16,5){$\bullet$}
\put(-4,5){$\bullet$}
\put(6,25){$\bullet$}
\put(8,7){\vector(-1,0){2}}
\put(3,17){\vector(1,2){2}}
\put(13,17){\vector(1,-2){2}}
\put(18,7){\line(-1,0){20}}
\put(18,7){\line(-1,2){10}}
\put(8,27){\line(-1,-2){10}}
\end{picture}
}
\end{picture}
&

\begin{picture}(160,50)
\put(0,25){
$) \;u(b)\,v(c)\,w(d)\;$}
\end{picture}
&
\begin{picture}(40,50)
\put(-95,25){$=(u\star v)\star w(a),$}
\end{picture}
\end{array}
\end{array}$$
the $\mu$-integration being taken over variables $b,c,d,t$ and $\tau$.
This leads, for $K=e^{iS}$, to an equality between two ``distribution valued  
functions'' on $M\times M\times M\times M$ (cf. formula~(\ref{ASSK})): for every 
quadruple of points $a,b,c,d$ in $M$ associativity for $\star$ reads
\begin{equation}\label{ASS}
\begin{array}{ccccc}
\begin{picture}(50,100)
\put(0,50){${\displaystyle \int}_M \exp\;i($}
\end{picture} & \carrehoriz{$t$} & 
\begin{picture}(80,100)
\put(0,50){$)\;\mu(t)={\displaystyle \int}_M \exp\;i($}
\end{picture} & \carrevertic{$\tau$} & 
\begin{picture}(50,100)
\put(0,50){$)\;\mu(\tau).$}
\end{picture}
\end{array}
\end{equation}
In the above formula, the diagram in the argument of the exponential in 
the LHS (respectively the RHS) stands for $S(a,b,t)+S(t,c,d)\;$ (respectively 
$S(a,d,\tau)+S(\tau,b,c)$).  

\begin{prop}\label{GEOM}
Let $(M,\mu,S)$ be a Weyl triple. Then, the associated three-point kernel 
$K=e^{iS}$ is geometrically associative.  
\end{prop}
\Pf
Fix four points $a,b,c,d$. Regarding \dref{coucou} and 
formula~(\ref{ASS}), one 
needs to construct our volume preserving diffeomorphism 
$\varphi:(M,\mu)\to(M,\mu)$ in such a way that for all $t$,
$$
\begin{array}{ccc}
\carrehoriz{$t$} & 
\begin{picture}(20,100)
\put(0,50){=}
\end{picture} & \carrevertic{$\varphi(t)$} 
\end{array}.
$$
We first observe 
that the data of four points $a,b,c,d$ determines what we call
an ``$S$-barycenter", that is a point $g=g(a,b,c,d)$ such that

$$
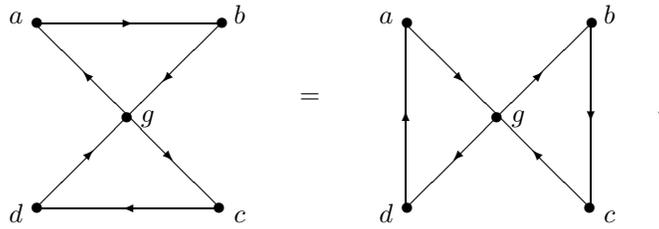
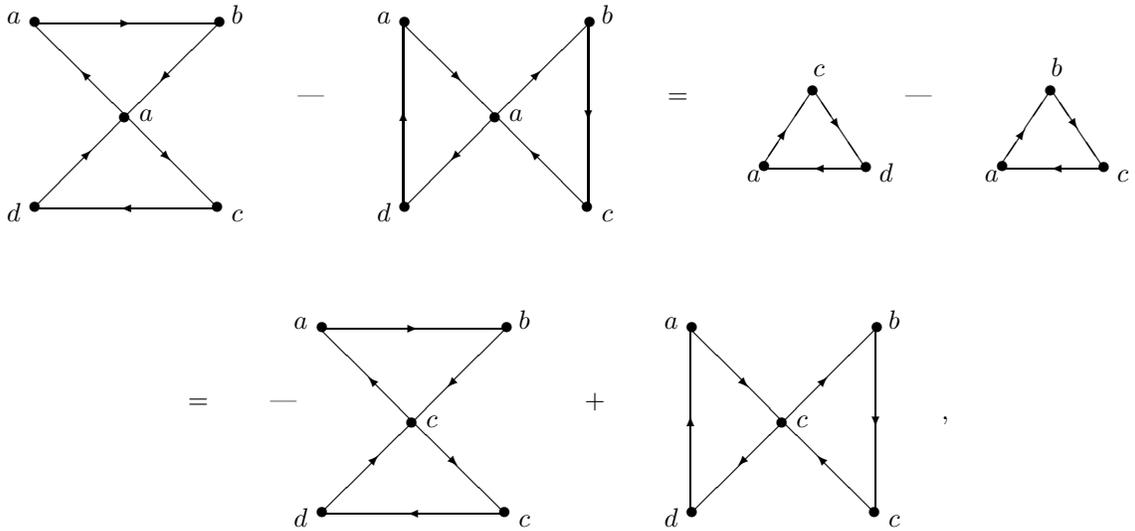
\begin{array}{ccc}
\carrehoriz{$g$} & 
\begin{picture}(20,100)
\put(0,50){=}
\end{picture} & \carrevertic{$g$} 
\end{array}.
$$
Indeed, since
$$
\begin{array}{ccccccc}
\carrehoriz{$a$} &
\begin{picture}(20,100)
\put(0,50){---}
\end{picture} & \carrevertic{$a$} &
\begin{picture}(10,100)
\put(0,50){=}
\end{picture} & \Triangle{a}{c}{d} & \begin{picture}(10,100)
\put(0,50){---}
\end{picture} & \Triangle{a}{b}{c}
\end{array}
$$

$$
\begin{array}{ccccc}
\begin{picture}(30,100)
\put(0,50){= \qquad ---}
\end{picture} &
\carrehoriz{$c$} &
\begin{picture}(20,100)
\put(0,50){+}
\end{picture} & \carrevertic{$c$}
\end{array},
$$
any continuous path joining $a$ to $c$ contains such a point $g$.

Now, we fix once for all such an $S$-barycenter $g$ for $\{ a,b,c,d\}$ and 
we adopt the following notation. For all $x$ and $y$ in $M$, the value 
of $S(g,x,y)$ is denoted by a ``thickened arrow''~:
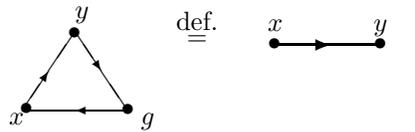
$$
\begin{array}{ccc}
\Triangle{x}{y}{g} & 
\begin{picture}(30,50)
\put(0,50){\mbox{   $\stackrel{\mbox{def.}}{=}$}}
\end{picture} & 
\begin{picture}(38,55)
\thicklines
\put(-2,55){$x$}
\put(38,55){$y$}
\put(-2,48){$\bullet$}
\put(38,48){$\bullet$}
\put(20,50){\vector(1,0){2}}
\put(0,50){\line(1,0){40}}
\end{picture}
\end{array}\quad.
$$
Again, a change of orientation in such an arrow changes the sign of its 
value. Also, property (iii) of a Weyl triple (Definition~\ref{WT}) which 
reads  
$$
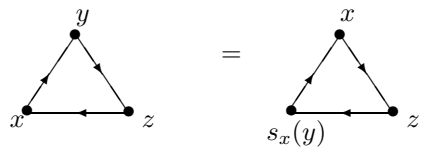
\begin{array}{ccc}
\begin{picture}(60,50)
\put(0,0){$x$}
\put(50,0){$z$}
\put(25,40){$y$}
\put(5,5){\trianglehorizontal}
\end{picture} & 
\begin{picture}(20,50)
\put(10,25){=}
\end{picture} & 
\begin{picture}(60,50)
\put(-3,-5){$s_x(y)$}
\put(50,0){$z$}
\put(25,40){$x$}
\put(5,5){\trianglehorizontal}
\end{picture}
\end{array},
$$
implies
$$
\begin{array}{ccc}
\begin{picture}(50,100)
\thicklines
\put(-2,55){$x$}
\put(38,55){$y$}
\put(-2,48){$\bullet$}
\put(38,48){$\bullet$}
\put(20,50){\vector(1,0){2}}
\put(0,50){\line(1,0){40}}
\end{picture}
&
\begin{picture}(10,90)
\put(0,50){=}
\end{picture} 
&
\begin{picture}(50,100)
\thicklines
\put(-2,55){$x$}
\put(38,55){$s_g(y)$}
\put(-2,48){$\bullet$}
\put(38,48){$\bullet$}
\put(20,50){\vector(-1,0){2}}
\put(0,50){\line(1,0){40}}
\end{picture}
\end{array}
$$
for all $x$ and $y$ in $M$.
While, from property (ii) (Definition~\ref{WT} with $m=g$), one gets
$$
\begin{array}{ccc}
\Triangle{x}{y}{z} &
\begin{picture}(10,90)
\put(0,50){=}
\end{picture} &
\Trianglefat{x}{y}{z}
\end{array}.
$$
Moreover, the barycentric property of $g$ reads
$$
\begin{array}{ccc}
\begin{picture}(100,100)
\thicklines
\put(10,10){\line(0,1){70}}
\put(80,10){\line(0,1){70}}
\put(10,45){\vector(0,1){2}}
\put(80,45){\vector(0,-1){2}}
\put(8,8){$\bullet$}
\put(8,78){$\bullet$}
\put(77,8){$\bullet$}
\put(78,78){$\bullet$}
\put(0,5){$d$}
\put(85,5){$c$}
\put(85,80){$b$}
\put(0,80){$a$}
\end{picture}
&
\begin{picture}(10,90)
\put(0,50){=}
\end{picture} 
&
\begin{picture}(100,100)
\thicklines
\put(10,10){\line(1,0){70}}
\put(10,80){\line(1,0){70}}
\put(45,10){\vector(-1,0){2}}
\put(45,80){\vector(1,0){2}}
\put(8,8){$\bullet$}
\put(8,78){$\bullet$}
\put(77,8){$\bullet$}
\put(78,78){$\bullet$}
\put(0,5){$d$}
\put(85,5){$c$}
\put(85,80){$b$}
\put(0,80){$a$}
\end{picture}.
\end{array}
$$
Hence
$$
\begin{array}{ccccc}
\carrehoriz{$t$}
&
\begin{picture}(10,90)
\put(0,50){=}
\end{picture} 
&
\begin{picture}(100,100)
\thicklines
\put(10,10){\line(1,0){70}}
\put(10,80){\line(1,0){70}}
\put(10,10){\line(1,1){70}}
\put(10,80){\line(1,-1){70}}
\put(30,30){\vector(1,1){2}}
\put(60,60){\vector(-1,-1){2}}
\put(45,10){\vector(-1,0){2}}
\put(45,80){\vector(1,0){2}}
\put(60,30){\vector(1,-1){2}}
\put(30,60){\vector(-1,1){2}}
\put(8,8){$\bullet$}
\put(8,78){$\bullet$}
\put(77,8){$\bullet$}
\put(78,78){$\bullet$}
\put(42,42){$\bullet$}
\put(0,5){$d$}
\put(85,5){$c$}
\put(85,80){$b$}
\put(0,80){$a$}
\put(50,43){$t$}
\end{picture}
&
\begin{picture}(10,90)
\put(0,50){=}
\end{picture} 
&
\begin{picture}(100,100)
\thicklines
\put(10,10){\line(0,1){70}}
\put(80,10){\line(0,1){70}}
\put(10,10){\line(1,1){70}}
\put(10,80){\line(1,-1){70}}
\put(30,30){\vector(1,1){2}}
\put(60,60){\vector(-1,-1){2}}
\put(10,45){\vector(0,1){2}}
\put(80,45){\vector(0,-1){2}}
\put(60,30){\vector(1,-1){2}}
\put(30,60){\vector(-1,1){2}}
\put(8,8){$\bullet$}
\put(8,78){$\bullet$}
\put(77,8){$\bullet$}
\put(78,78){$\bullet$}
\put(42,42){$\bullet$}
\put(0,5){$d$}
\put(85,5){$c$}
\put(85,80){$b$}
\put(0,80){$a$}
\put(50,43){$t$}
\end{picture}
\end{array}
$$

$$
\begin{array}{cccc}
\begin{picture}(10,90)
\put(0,50){=}
\end{picture} 
&
\carreverticfat{$s_g(t)$}
&
\begin{picture}(10,90)
\put(0,50){=}
\end{picture} 
&
\carrevertic{$s_g(t)$}.
\end{array}
$$
One can therefore choose our diffeomorphism $\varphi$ as 
$$
\varphi=s_g.
$$
\EPf

\begin{rmks}\label{rmk:Ravel}\hspace{-.1cm}
{\rm \begin{enumerate}

\item[(i)] Given four points $\{a,b,c,d\}$ in the Euclidean plane $E^{2}$, 
one can define two smooth functions $F$ and $G\in C^{\infty}(E^{2},\R)$ 
as follows. For $t\in E^{2}$, $G(t)$ is the sum of the signed Euclidean 
areas of the oriented Euclidean triangles ${\Delta abt}$ and ${\Delta dtc}$, 
i.e.  
$$
G(t)=SA({\Delta abt})+SA({\Delta dtc})
$$
where, for a closed path $\gamma$ in $E^{2}$, one sets $SA (\gamma) = 
\int_{\gamma}x\, dy$.\\    
Similarly, one defines for the other pair of triangles ${\Delta dat}$ 
and ${\Delta bct}$,
$$
F(t)=SA({\Delta dat})+SA({\Delta bct})
$$
or equivalently, if $SA(abcd)$ denotes the signed area of the oriented 
quadrilateron with vertices $a,b,c,d$, 
$$
F(t)=SA(abcd)-G(t).
$$
If our four points are the vertices of a {\sl square} in $E^{2}$, 
one obviously has the equality
$$
F=G.
$$
For arbitrary quadrilaterons, equality between $F$ and $G$ does not 
generally hold. But the proof of Proposition~\ref{GEOM} leads to a natural 
generalization~:  

\vspace{1mm}

{\sl For every quadruple of points $\{a,b,c,d\}$ in $E^{2}$, one has 
$$
F=s_{g}^{\star}G
$$
where $g$ is the center of mass of the quadrilateron $abcd$, and, 
where $s_g:E^2\to E^2$ is the Euclidean symmetry $s_g(x)=2g-x$.}

\vspace{1mm}

Proposition~\ref{GEOM} tells us that this elementary property of Euclidean quadrilaterons of 
$E^{2}$ {\bf is} the geometric content of associativity of Weyl's product. 
Indeed, associativity for the Weyl product (on $\R^2$) reads 
(cf.~formula~(\ref{ASS}))~:   
$$
\int_{\R^2}\exp\, iF(t)dt=\int_{\R^2}\exp\, iG(\tau)d\tau.
$$
This equality being realized by the change of variables $\tau = 
s_{g}(t)$.   

\item[(ii)]
Proposition~\ref{GEOM} also leads to a geometric understanding of the 
associativity of Rieffel's deformed product obtained from an action 
of $\R^d$ on a $C^*$-algebra \cite{R1}. This is based on the following 
observation.

Let $G$ be a Lie group, endowed with a left invariant Haar measure $\mu$. 
Assume that there exists a geometrically associative three-point kernel 
$K \in C^{\infty}(G \times G \times G)$ (cf.~\dref{coucou}) on $(G,\mu)$ 
which is invariant under the diagonal action of $G$ on $G \times G \times 
G$ by left translations. Then, the data of an action $\alpha : G \times A 
\to A$ of $G$ on some (topological) associative algebra $A$ by 
automorphisms yields a new associative product on $A$, provided the situation 
is sufficiently regular. Indeed, on just observes that the geometric 
associativity of $K$ implies, at least formally, associativity of the product 
$\star$ on $A$ defined by the formula~: 
\begin{equation}\label{formule}
a \star b = \int_{G \times G} K(e,g,h) \alpha_g(a) \alpha_h(b),
\end{equation}
where $a, b \in A$, and where integration is taken over $g$ and $h$
with respect to the 
measure $\mu \otimes \mu$ on $G \times G$. In the case of an action of 
the Abelian group $G = \R^d$, endowed with the Euclidean scalar 
product $<,>$, every skewsymmetric matrix $J \in so(\R^d)$ yields an 
invariant three-point function $S^J \in C^\infty(\R^d \times \R^d \times 
\R^d)$ via the formula 
\begin{equation}
S^J(x,y,z) = <x,Jy> + <y,Jz> + <z,Jx>.
\end{equation}
Defining the Euclidean symmetries of $\R^d$ as $s_x(y) = 2 x - y$, 
one gets a Weyl triple $(\R^d, {\rm Haar}, S^J)$ (\dref{WT}). Hence, 
by \pref{GEOM}, an invariant geometrically associative kernel 
$$
K = e^{i S^J}.
$$
When given an action $\alpha : \R^d \times A \to A$, formula (\ref{formule}) 
with $K = e^{i S^J}$ yields Rieffel's product in \cite{R1}.

\item [(iii)] In the case of the hyperbolic plane endowed with its natural structure 
of symmetric space, it is tempting to ask whether the three-point function 
defined by the symplectic area $SA(x,y,z)$ of the geodesic triangle ${\Delta 
xyz}$ would satisfy properties (i)---(iii) in \dref{WT}. The answer is 
negative. Indeed, properties (ii) and (iii) would imply unboundedness of 
$SA$. This indicates that, when requiring some compatibility between $S$ 
and the symmetries, property (iii) is somehow too strong.
\end{enumerate}
}\end{rmks}

Proposition~\ref{GEOM} and Remark~\ref{rmk:Ravel} lead us to consider a class of manifolds which carry a 
large family of ``symmetries''. Symplectic symmetric spaces are 
defined in \sref{SSS}. They constitute a class of affine 
symplectic manifolds. We will now make precise  what we mean by
``WKB-quantization'' in the context of affine 
symplectic manifolds.

Let $(M,\omega,\nabla)$ be a $2n$-dimensional {\bf affine symplectic 
manifold}, that is, $(M, \omega)$ is a smooth connected symplectic manifold 
and $\nabla$ is an torsion-free affine  connection on $M$ such that 
$\nabla\omega=0$. Its {\bf automorphism group} $\mbox{Aut}(M,\omega,\nabla)$ 
is defined as  
$$
\mbox{Aut}(M,\omega,\nabla)=\mbox{Aff}(\nabla)\cap\mbox{Symp}(\omega)
$$
where $\mbox{Aff}(\nabla)$ is the group of affine transformations of 
the affine manifold $(M,\nabla)$ and where $\mbox{Symp}(\omega) $
denotes the group of symplectomorphisms of $(M,\omega)$. Note that, 
since $\mbox{Aff}(\nabla)$ is a Lie group of transformations of $M$ 
(cf.~\cite{KN}), so is $\mbox{Aut}(M,\omega,\nabla)$.\\

\begin{dfn}\label{WKB}
Let $G$ be a subgroup of $Aut(M,\omega,\nabla)$. A {\bf $G$-invariant  
WKB-quantization} of $(M,\omega,\nabla)$ is a triple 
$(\{(\e_\hbar,\star_\hbar)\}_{\hbar\geq 0}, S, \{a_\hbar\}_{\hbar\geq 0})$ satisfying the following 
properties.
\begin{enumerate}
\item[(i)] $\{(\e_\hbar,\star_\hbar)\}_{\hbar\geq 0}$ is a one-parameter family of 
	associative ${}^*$-algebras such that~:
	\begin{enumerate}
	\item[(i.1)] $(\e_0,\star_0)$ is a Poisson subalgebra of 	
	$C^\infty(M)$ endowed with the usual pointwize multiplication
	of functions (the Poisson structure $\{\, , \,\} $ is the one
	associated to the symplectic form $\omega$).
	\item[(i.2)] For all $\hbar\geq 0$, $\e_{\hbar}$ is a ${}^*$-linear 
	subspace of $C^\infty(M)$ such that the following inclusions hold
	$$ \d(M)\subset \e_0 \subset \e_{\hbar},$$ where 
	$\d(M)$ denotes the space of smooth compactly supported functions on $M$ 
	and where on $C^{\infty}(M)$ 
	the involution is the complex conjugation.  
	\end{enumerate}
\item[(ii)] $S$ is a real valued smooth three-point function~: $S\in 
C^\infty(M \times M \times M, \R)$ such that, for all $x_0\in M$ the partial 
function $S(x_0,\, . \,,\, .  \,)\in C^\infty(M \times M, \R)$ has an 
nondegenerate critical point at $(x_0,x_0)\in M\times M$. One furthermore 
requires the function $S$ to be invariant under the diagonal action of 
$G$ on $M\times M \times M$. 
\item[(iii)] $\{a_\hbar\}_{\hbar\geq 0}$ is a smooth (with respect to
$\hbar$) family of positively real valued three-point 
functions~:
$a_\hbar\in C^\infty(M\times M \times M, \R^+)$ which are invariant
under the diagonal action of 
$G$ on $M\times M \times M$.
\item [(iv)] At the level of the subspace $\d(M)\subset \e_\hbar$ 
($\hbar>0$),  the multiplication $\star_\hbar$
reads 
$$
(u\star_\hbar v)(x)=\frac{1}{\hbar^{2n}}\int_{M\times M}a_\hbar(x,y,z)\, \exp\left( 
\frac{i}{\hbar}\,S(x,y,z)\right)\, u(y)\, v(z) \, dy\, dz
$$
where $u,v\in \d(M)$ and where $dy$ (and $dz$) stands for the Liouville measure 
$\frac{\omega^n}{n!}$.
\item[(v)] For all $x\in M$ and all	$u,v\in \d(M)\subset \e_\hbar$ 
supported in a sufficiently small neighborhood of $x$, a stationary phase 
method yields the following asymptotic expansion (cf.~\cite{FM})~:
$$
(u\star_\hbar v)(x)\sim u(x)v(x) + \frac{\hbar}{i}c_1(u,v)(x) + o(\hbar^2)
$$ 
with 
$$
\frac{1}{2}\left(c_1(u,v)(x)-c_1(v,u)(x)\right)=\{u,v\}(x).
$$
\end{enumerate}
\end{dfn}
This definition is very much inspired from Rieffel's definition of strict 
quantization. However, the topological framework is weaker than Rieffel's 
one. For instance, the deformed algebras are function algebras but they 
do not a priori carry any topological structure. Also Dirac's condition 
does only hold at the formal level (see item $(v)$ in \dref{WKB}). 
Note moreover that each function space $\e_{\hbar}$ is not required 
to be invariant under the action of $G$ on $C^\infty(M)$. However, the 
space of compactly supported functions $\d(M)$ is $G$-invariant, and, by
invariance of the functions $a_{\hbar}$ and $S$, one will always have
$$
g(u\star_{\hbar}v)=(gu)\star_{\hbar}(gv)\quad\forall g\in G
$$
as soon as $u,v\in\d(M)\subset \e_{\hbar}$.

\section{General facts about symplectic symmetric spaces}\label{SSS}

\begin{dfn}\label{SSSDEF}\cite{Bi3, Bi2} A {\bf symplectic symmetric space} 
is a triple $(M, \omega,s)$, where $(M,\omega)$ is a smooth connected 
symplectic manifold, and where $s : M \times M \to M$ is a smooth map such 
that 
  
\begin{enumerate} 
\item[(i)] for all $x$ in $M$, the partial map $s_x : M \to M : y \mapsto 
s_x (y) := s(x,y)$ is an involutive symplectic diffeomorphism of $(M,\omega)$ 
called the {\bf symmetry} at $x$. 
\item[(ii)] For all $x$ in $M$, $x$ is an isolated fixed point of $s_x$. 
\item[(iii)] For all $x$ and $y$ in $M$, one has $s_xs_ys_x=s_{s_x (y)}$. 
\end{enumerate}

\end{dfn}

\begin{dfn} Two symplectic symmetric spaces  $(M,\omega,s)$ 
and $(M',\omega ',s')$ are {\bf isomorphic} if there exists a symplectic 
diffeomorphism $\varphi: (M,\omega) \rightarrow (M',\omega')$ such that 
$\varphi s_x=s'_{\varphi (x)} \varphi$. Such a $\varphi $ is called an {\bf
isomorphism} of $(M,\omega,s)$ onto $(M',\omega',s')$. When 
$(M,\omega,s) = (M',\omega',s')$, one talks about {\bf automorphisms}. 
The group of all automorphisms of the symplectic symmetric space 
$(M,\omega,s)$ is denoted by $Aut(M,\omega,s)$.
\end{dfn}

\begin{prop} On a symplectic symmetric space $(M,\omega,s)$, there exists 
one and only one affine connection $\nabla$ which is invariant under the 
symmetries. Moreover, this connection satisfies the following properties. 

\begin{enumerate} 
\item[(i)] For all smooth tangent vector fields $X,Y,Z$  on $M$ and all 
points $x$ in $M$, one has $$ \omega_{x} (\nabla _X Y,Z)=\frac{1}{2}
X_x.\omega (Y+s_{x_{\star}}Y,Z).$$ \item[(ii)] $(M,\nabla )$ is an affine 
symmetric space. In particular $\nabla$ is torsion free and its curvature 
tensor is parallel.  
\item[(iii)] The symplectic form $\omega$ is parallel; $\nabla$ is therefore 
a symplectic connection.   
\item[(iv)] One has
$$\mbox{Aut}(M,\omega,s)=\mbox{Aut}(M,\omega,\nabla)=\mbox{Aff}(\nabla)\cap\mbox{Symp}(\omega).$$
\end {enumerate} 

\end{prop}

The connection $\nabla$ on the symmetric space $(M,s)$ is called the {\bf 
Loos connection}. The following facts are classical (see \cite{Lo}, 
v.~I, \cite{KN}, v.~II, Chapters X and XI).  

\begin{thm}\label{TG}
Let $(M,\omega,s)$ be a symplectic symmetric space and $\nabla$ its 
Loos connection. Fix $o$ in $M$ and denote by $H$ the stabilizer of 
$o$ in $Aut(M,\omega,s)$. Denote by $G$ the {\bf transvection group} of 
$(M,s)$ (i.e.~the subgroup of $Aut(M,\omega,s)$ generated by $\{ s_x \circ 
s_y \, ; \, x,y \in M \}$)  and set $K=G \cap H$. Then, 

\begin{enumerate} 
\item[(i)] the transvection group $G$ turns out to be a connected Lie transformation 
group of $M$. It is the smallest subgroup of $Aut(M,\omega,s)$ which is transitive 
on $M$ and stabilized by the conjugation $\tilde{\sigma} : 
Aut(M,\omega,s) \to Aut(M,\omega,s)$ defined by $\tilde{\sigma}(g)=s_o g s_o$.
\item[(ii)] The homogeneous space $M=G/_{\textstyle{K}}$ is reductive. The 
Loos connection $\nabla$ coincide with the canonical connection induced 
by the structure of reductive homogeneous space. 
\item[(iii)] Denoting by $G^{\tilde{\sigma}}$ the set of $\tilde{\sigma}$-fixed 
points in $G$ and by $G_0^{\tilde{\sigma}}$ its neutral connected component, 
one has $$ G_0^{\tilde{\sigma}} \subset K \subset G^{\tilde{\sigma}}.$$ 
Therefore, the Lie algebra of $K$ is $\k$. Moreover, it is isomorphic to the holonomy algebra 
with respect to the canonical connection $\nabla$.
\item[(iv)] Denote by $\sigma$ the involutive automorphism of the Lie 
algebra $\g$ of $G$ induced by the automorphism $\tilde{\sigma}$. Denote 
by $\g = \k \oplus \p$ the decomposition in $\pm 1$-eigenspaces for $\sigma$. 
Then, identifying $\p$ with $T_o(M)$, one has $$ exp (X) = 
s_{Exp_o(\frac{1}{2} X)} \circ s_o $$ for all $X$ in a neighborhood of 
$0$ in $\p$. Here $exp$ is the exponential map $exp : \g \to G$ and 
$Exp_o$ is the exponential map at point $o$ with respect to the connection 
$\nabla$.   

\end{enumerate}
\end{thm}

\begin{dfn} Let $(\g,\sigma)$ be an {\bf involutive algebra}, that is, 
${\cal G}$ is a finite dimensional real Lie algebra and $\sigma$ is an 
involutive automorphism of ${\cal G}$. Let $\Omega$ be a skewsymmetric 
bilinear form on $\g$. Then the triple $({\cal G},\sigma, \Omega)$ is called 
a {\bf symplectic triple} if the following properties are satisfied. 

\begin{enumerate}	
\item[(i)]	Let ${\cal G}={\cal K}\oplus {\cal P}$   where ${\cal K}$ (resp. 
${\cal P}$) is the $+1$ (resp. $-1$) eigenspace of $\sigma$. Then $[{\cal	
P},{\cal P}]={\cal K}$ and the representation of ${\cal K}$ on ${\cal P}$, 
given by the adjoint action, is faithful. 
\item[(ii)] $\Omega$ is	a Chevalley 2-cocycle for the trivial representation 
of ${\cal G}$ on $\R$ such that for any $X$ in ${\cal K}$, 
$i(X){\Omega}=0$. Moreover, the restriction of $\Omega$ to ${\cal P} \times 
{\cal P}$ is nondegenerate.
\end{enumerate} 

The dimension of ${\cal P}$ defines the {\bf dimension}	of	the triple. 
Two	such triples $({\cal G}_i,\sigma_i,{\Omega}_i)$ $(i=1,2)$ are 
{\bf isomorphic} if there exists a Lie algebra isomorphism 
$\psi :{\cal G}_1\rightarrow{\cal G}_2$ such that $\psi \circ \sigma_1 = 
\sigma_2 \circ \psi$ and $\psi^*{\Omega}_2={\Omega}_1$. 
\end{dfn}

\tref{TG} associates to a symplectic symmetric space $(M,\omega,s)$ 
an involutive Lie algebra $(\g,\sigma)$. Denoting by $\pi : G \to M$ 
the natural projection, one checks that the triple $(\g, \sigma, \Omega = 
\pi^* (\omega_o))$ is a symplectic triple. This implies the next 
proposition.  

\begin{prop} There is a bijective correspondence between the isomorphism
classes of simply connected symplectic symmetric spaces $(M,\omega,s)$ 
and the isomorphism classes of symmetric triples $({\cal G},\sigma,{ 
\Omega})$.  
\end{prop}

Since a symmetric symplectic manifold $(M,\omega,s)$ is a symplectic 
homogeneous space of its transvection group $G$, it seems natural, when 
possible, to relate $(M,\omega,s)$ to a coadjoint orbit of $G$ in ${\cal 
G}^\star$. Recall first the two following definitions.  

\begin{dfn} Let $G$ be a Lie group of symplectomorphisms acting on a 
symplectic manifold $(M,\omega)$. For every element $X$ in the Lie algebra 
$\g$ of $G$, one denotes by $X^*$ the {\bf fundamental vector field} 
associated to $X$, i.e.~for $x$ in $M$,
$$
X^*_x = \frac{d}{dt} \exp (-tX) x|_{t=0}.
$$
The action is called {\bf weakly Hamiltonian} if for all $X$ in $\g$ 
there exists a smooth function $\lambda_X \in C^\infty(M)$ such that 
$$
i(X^*)\omega = d \lambda_X.
$$ 
In this case, if the correspondence $\g \to C^\infty (M) : X \mapsto 
\lambda_X$ is also a homomorphism of Lie algebras, one says that the action 
of $G$ on $(M,\omega)$ is {\bf Hamiltonian}. (The Lie algebra structure 
on $C^\infty (M)$ is defined by the Poisson bracket.)
\end{dfn}

\begin{prop} Let $t=({\cal G},\sigma,\Omega)$ be a symplectic triple 
and let $(M,\omega, s)$ be the associated simply connected symplectic 
symmetric space. The action of the transvection group $G$ on $M$ is 
Hamiltonian if and only if ${\Omega}$ is a Chevalley coboundary, that 
is, there exists an element $\xi$ in $\g^*$ such that $\Omega = 
\delta \xi$. In this case, $(M,\omega,s)$ is a $G$-equivariant symplectic 
covering of ${\cal O}$, the coadjoint orbit of $\xi$ in ${\cal G}^\star$. 
\end{prop}

The action of the transvection group $G$ is in general not Hamiltonian. 
We therefore need to consider a one-dimensional central extension of $G$ 
rather than $G$ itself. At the infinitesimal level, this corresponds to 
extending the algebra $\g$ by the 2-cocycle $\Omega$. This way, one associates 
to any symplectic symmetric space an exact triple in the following sense 
(see \cite{Bi1} for details).

\begin{dfn}\label{ET} An {\bf exact triple} is a triple $\tau = 
(\g,\sigma,{\Omega})$, where 

\begin{enumerate} 
\item[(i)] $(\g,\sigma)$ is an involutive Lie algebra such that, if 
$\g=\KK \oplus \p$ is the decomposition with respect to $\sigma$, one has   
$\left[\p,\p\right]=\KK$. 
\item[(ii)] ${\Omega}$ is a Chevalley 2-coboundary (i.e. ${\Omega}=\delta 
\xi \,, \quad \xi \in \g^\star$) such that $i(\KK){\Omega}=0$ and 
${\Omega}|_{\p \times \p}$ is symplectic. 
\end{enumerate} 

\end{dfn} 

\begin{rmk}\label{DIM(Z)}{\rm 
\begin{enumerate}
\item[(i)] One can choose $\xi \in \g^\star$ such that $\xi(\p)=0$. One 
therefore writes, with a slight abuse of notation, $\xi\in\KK^\star$.
\item[(ii)] Observe that, when associated to a (transvection) symplectic 
triple, the center $\z(\g)$ of the Lie algebra $\g$ occurring in an exact 
triple is at most one dimensional. Indeed, on the one hand, exactness 
implies $\z(\g)\subset\k$. One the other hand, faithfulness of the holonomy 
representation forces $\dim(\z(\g)\cap\k)\leq1$ since $\k$ is either the 
holonomy algebra itself or a one dimensional central extension.
\end{enumerate}
}\end{rmk}

\section{Phase functions on symplectic symmetric spaces}\label{PHSSS}

Motivated by the definition of a Weyl triple (\dref{WT}), as well as by 
the third part of \rsref{rmk:Ravel}, we now make the following definition, 
adapting to symmetric spaces the notion of phase function.

\begin{dfn}\label{PHASE}
Let $(M,\omega,s)$ be a symplectic symmetric space (see~\dref{SSSDEF}). 
A smooth function $S:M\times M\times M \to \R$ satisfying the 
following properties~:
\begin{enumerate}
\item[(i)] $S(x,y,z)=S(z,x,y)=-S(y,x,z)$;
\item[(ii)] the function $S$ is invariant under the symmetries i.e. 
$$
S(s_m(x),s_m(y),s_m(z))=S(x,y,z)\qquad \forall x,y,z,m\in M;
$$
\item[(iii)] for all $x \in M$, the symmetry $s_x:M\to M$ is such that
$$
S(x,y,z)=-S(x,s_x(y),z)\qquad \forall y,z\in M,
$$
\end{enumerate}
is called an {\bf admissible} function.
\end{dfn}

In \cite{W1}, Weinstein proved that, in the case of a Hermitian symmetric 
space $(M,\omega,s)$ of the noncompact type, the phase function $S_W$, 
occurring in the expression of a given invariant WKB-quantization of 
$(M,\omega)$ defined by an oscillatory integral formula of the type
\begin{equation}\label{SW}
u\star v(x)=\int u(y)v(z) a_\hbar(x,y,z) 
e^{\frac{i}{\hbar}S_W(x,y,z)}dy\,dz,
\end{equation}
must be as follows. 

\begin{enumerate}
\item[$\bullet$]
Let $x,y$ and $z$ be three points in $M$ such that 
the following equation admits a solution $t$~:  
$$
t=s_xs_ys_z(t)
$$
($t$ is unique if it exists). 
\item[$\bullet$] 
Let $\Sigma$ be a surface in $M$ 
bounded by the geodesic triangle ${\Delta tAB}$ 
where $A=s_x(t)$ and $B=s_y(A)=s_ys_x(t)\quad$ ($t=s_z(B)$). 
\end{enumerate}

Then, if for 
some formal amplitude of the form 
$$
a_\hbar(x,y,z)=a_0(x,y,z) +\hbar a_1(x,y,z) + \hbar^2 a_2(x,y,z) + 
...\quad,
$$
the product (\ref{SW}) defines an invariant deformation quantization of 
$(M,\omega)$, the value of the ``WKB"-phase function 
$S_W$ on $(x,y,z)$ is given by 

$$
S_W(x,y,z)=-\int_\Sigma\omega.
$$

Practically, the function $S_W$ is hard to compute explicitly; however, 
some two dimensional examples have been treated in \cite{Q}. The problem 
of finding the amplitude $a_\hbar$ is open.   

\begin{prop}
The function $S_W$ is admissible.
\end{prop}
\Pf Let $x,y$ and $z$ be three points in $M$, and set
$$
\left\{
\begin{array}{ccl}
Y & = & s_zs_ys_x(Y)\qquad(*)\\
Z & = & s_x(Y)\\
X & = & s_ys_x(Y)=s_y(Z)
\end{array}
\right.
$$
and 
$$
\left\{
\begin{array}{ccl}
\zeta & = & s_{s_x(y)}s_zs_x(\zeta)\qquad(**)\\
\xi & = & s_x(\zeta)\\
\eta & = & s_zs_x(\zeta)=s_z(\xi),
\end{array}
\right.
$$
(provided equations (*) and (**) admit solutions). The only thing we 
really need to show is that the function $S_{W}$ satisfies property 
(iii) in \dref{PHASE}; that is,
$$
SA(X,Y,Z) = SA(\xi, \eta, \zeta),
$$
where the function $SA$ is defined as follows. For a sequence of points $\{x_i\}_{0\leq i\leq N}$, the 
expression $SA(x_0,...,x_N)$ means $SA(\gamma) = \int_{\gamma}\alpha$, 
where $\gamma$ is the piecewise geodesic  path whose $i^{\mbox{th}}$ geodesic 
segment starts from point $x_i$ and ends at point $x_{i+1\mbox{mod}(N+1)}$, 
and where $\alpha$ is a $1$-form such that $\omega = d\alpha$.

On the first hand, one has~: 
$
s_xs_{s_x(y)}s_z(X)=s_xs_xs_ys_xs_z(X)=s_ys_xs_z(X)=X.
$
Hence $X=\xi$ and $s_{s_x(z)}(Z)=\zeta$. On the second hand, the 
invariance of $SA$ under the symmetries yields~:
$SA(\xi,\eta,\zeta)=SA(s_x(\xi),s_x(\eta),s_x(\zeta))=SA(\zeta,Z,\xi)$. 
Moreover, $SA(\xi,\eta,\zeta)+SA(\zeta,Z,\xi) = SA(\xi,\eta,\zeta,Z)$, 
hence 
$$S_W(x,s_x(y),z)=\frac{1}{2}SA(\xi,\eta,\zeta,Z).$$ Similarly, one 
gets $$S_W(x,y,z)=-\frac{1}{2}SA(\xi,\eta,\zeta,Z).
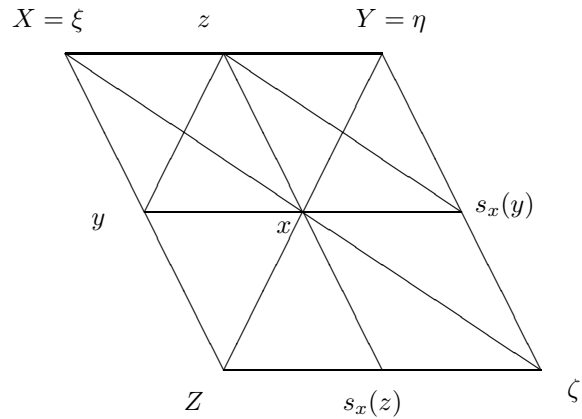
$$

$$
\begin{picture}(220,180)
\put(20,140){\line(1,-2){60}}
\put(80,140){\line(1,-2){60}}
\put(140,140){\line(1,-2){60}}
\put(20,140){\line(1,0){120}}
\put(80,20){\line(1,0){120}}
\put(50,80){\line(1,0){120}}
\put(20,140){\line(3,-2){180}}
\put(80,20){\line(1,2){60}}
\put(50,80){\line(1,2){30}}
\put(80,140){\line(3,-2){90}}
\put(0,150){$X=\xi$}
\put(70,150){$z$}
\put(130,150){$Y=\eta$}
\put(30,75){$y$}
\put(100,72){$x$}
\put(175,80){$s_x(y)$}
\put(65,5){$Z$}
\put(125,5){$s_x(z)$}
\put(210,10){$\zeta$}
\end{picture}
$$

\EPf

Now, suppose we are given, on a symplectic symmetric space $(M,\omega,s)$, an 
admissible three-point function $S$. The choice of a base point   
$o$ in $M$ then determines a two-point function $u$ on $M \times M$~:
$$
u(x,y)=S(o,x,y).
$$
The following proposition tells us, in the case where $(M,\nabla)$ is 
strictly geodesically convex, when a given two-point function comes 
from an admissible three-point function.  

\begin{prop}\label{MIDPOINT}
Let $(M,\omega,s)$ be a symplectic symmetric space and let $o$ be a 
point in $M$. Assume the existence of a smooth {\bf midpoint map}, that is, 
a map $M\to M: x\to\frac{x}{2}$ such that $s_{\frac{x}{2}}(o)=x$. Denote 
by ${\textstyle Stab}(o)$ the stabilizer of $o$ in the automorphism group 
of $(M,\omega,s)$. Let $u:M\times M\to \R$ be a ${\textstyle 
Stab}(o)$-invariant smooth two-point function such that
$$
u(x,y)=-u(y,x)=-u(x,s_x(y))=-u(x,s_{\frac{x}{2}}(y))\qquad\forall 
x,y\in M.
$$
Then, the three-point function $S$ defined by 
$$
S(x,y,z)=u(s_{\frac{x}{2}}(y),s_{\frac{x}{2}}(z))
$$
is admissible.
\end{prop}
\Pf
The function $S$ is invariant under the symmetries. Indeed, let $g$ be an 
element of the group of transformations of $M$ 
generated by
the symmetries. For all $x$ in $M$, one has $s_{\frac{gx}{2}}(o)=gx$ 
that is $g^{-1}s_{\frac{gx}{2}}(o)=s_{\frac{x}{2}}(o)$;
hence $s_{\frac{x}{2}}g^{-1}s_{\frac{gx}{2}}$ is an element of 
${\textstyle Stab}(o)$. By invariance under the action of
the stabilizer, one has 
$$u(s_{\frac{gx}{2}}(gy),s_{\frac{gx}{2}}(gz))=u(s_{\frac{x}{2}}g^{-1}s_{\frac{gx}{2}}
s_{\frac{gx}{2}}(gy), 
s_{\frac{x}{2}}g^{-1}s_{\frac{gx}{2}}s_{\frac{gx}{2}}(gz))=u(s_{\frac{x}{2}}(y),
s_{\frac{x}{2}}(z)).$$

Now, since 
$s_{\frac{x}{2}}s_{\frac{z}{2}}s_{\left(\frac{s_{\frac{z}{2}}(x)}{2}\right)}$ 
is an element of 
${\textstyle Stab}(o)$, 
one has
$$
u(s_{\frac{z}{2}}(x),s_{\frac{z}{2}}(y))=
-u(s_{\frac{z}{2}}(x),s_{\left(\frac{s_{\frac{z}{2}}(x)}{2}\right)}s_{\frac{z}{2}}(y))=
$$

$$
=-u(s_{\frac{x}{2}}s_{\frac{z}{2}}s_{\left(\frac{s_{\frac{z}{2}}(x)}{2}\right)}
s_{\frac{z}{2}}(x),
s_{\frac{x}{2}}s_{\frac{z}{2}}s_{\left(\frac{s_{\frac{z}{2}}(x)}{2}\right)}s_{\left(
\frac{s_{\frac{z}{2}}(x)}{2}\right)}
s_{\frac{z}{2}}(y))=
-u(s_{\frac{x}{2}}s_{\frac{z}{2}}s_{\left(\frac{s_{\frac{z}{2}}(x)}{2}\right)}
s_{\frac{z}{2}}
s_{\frac{x}{2}}(o),s_{\frac{x}{2}}(y)) =
$$

$$
=-u(s_{\frac{x}{2}}s_{\frac{z}{2}}\left(
s_{\frac{x}{2}}s_{\frac{z}{2}}s_{\left(\frac{s_{\frac{z}{2}}(x)}{2}\right)} 
\right)^{-1}(o),s_{\frac{x}{2}}(y))
=-u(s_{\frac{x}{2}}(z),s_{\frac{x}{2}}(y))=u(s_{\frac{x}{2}}(y),s_{\frac{x}{2}}(z)).
$$
This shows that the function $S$ satisfies the condition $(i)$ in 
\dref{PHASE}. The remaining condition is obviously satisfied. 

\EPf

\begin{dfn}\label{TPF} Let $(M,\omega,s)$ be a symplectic symmetric space 
admitting a midpoint map (cf.~\pref{MIDPOINT}) with respect to a 
point $o$ in $M$. An {\bf $o$-admissible two-point function} on $M$ is a 
${\textstyle Stab}(o)$-invariant smooth function $u:M\times M\to \R$  such 
that 
$$
u(x,y)=-u(y,x)=-u(x,s_x(y))=-u(x,s_{\frac{x}{2}}(y))\qquad\forall 
x,y\in M.
$$
\end{dfn}
It turns out that admissible two-point functions are easy to determine.

\begin{ex}\label{EX2DIM}
{\rm
The triple $(M,\omega,s)=(\R^2,dp\wedge dq,s)$ with 
$$
s_{(p,q)}(p',q') = (2p-p',2\cosh(p-p')q-q')
$$ 
defines a (non-metric) symplectic symmetric space \cite{Bi1}, \cite{Bi3}. 
It is strictly geodesically convex. One checks that the function $u:M\times 
M\to \R$ defined by 
$$
u((p,q),(p',q')) =  \sinh(p')q - \sinh(p)q'
$$ 
is an $(0,0)$-admissible two-point function whose corresponding admissible 
three-point function $S$ is given by~: 
$$
S((p,q),(p',q'),(p",q"))=\sinh(p-p')q"+\sinh(p"-p)q'+\sinh(p'-p")q.
$$
}\end{ex}

\section{Elementary solvable symplectic symmetric spaces}\label{ESSS}

In \dref{ESET} below, we define a particular type of solvable 
symmetric spaces that we call elementary. It has been proven 
(\cite{Bi1}, Proposition~3.2) that every solvable symmetric space can be 
expressed as the result of a sequence of split extensions by Abelian 
(flat) factors successively taken over an elementary solvable 
symmetric space. We therefore consider elementary solvable symmetric 
spaces as the ``first induction step'' when studying solvable symmetric 
spaces.

\begin{dfn}\label{ESET}
A symplectic symmetric space $(M,\omega,s)$ is called an {\bf elementary 
solvable} symplectic symmetric space if its associated exact triple 
$(\g,\sigma,\Omega = \delta \xi)$ (see~\dref{ET}) is of the following type.
\begin{enumerate}
\item[(i)] The Lie algebra $\g$ is a split extension of Abelian Lie 
algebras $\a$ and $\b$~:
$$
\b\longrightarrow\g\stackrel{\longleftarrow}{\longrightarrow}\a.
$$
\item[(ii)] The automorphism $\sigma$ preserves the splitting 
$\g=\b\oplus\a$.
\end{enumerate}
Such an exact triple (associated to an elementary solvable symplectic 
symmetric space) is called an {\bf elementary solvable exact triple} 
(briefly: {\bf ESET}). 
\end{dfn}

Observe that, since $\a \cap \KK \subset \a \cap [\g,\g] = 0$, one 
has $\a \subset \p$. Therefore $\b = \KK \oplus \EK$, with $\EK 
\subset \p$. Moreover, since $\EK$ and $\a$ are Abelian and 
$\Omega$ is nondegenerate, the subspaces $\a$ and $\EK$ of $\p$ are 
Lagrangians in duality.

Now, let $\rho:\a\to\mbox{End}(\b)$ be the splitting homomorphism. 
The group law on $G$ identified with $\g=\a\times\b$ is given by
$$
(a,b).(a',b')=(a+a',\exp(\rho(a))b'+b).
$$
The (simply connected) symplectic symmetric space $(M,\omega,s)$ 
associated to an ESET $(\g,\sigma,\Omega)$ is described as follows 
\cite{Bi3, Bi1}.
 
\begin{prop}\label{SPACES}
\begin{enumerate}
\item[(i)] The homogeneous space $G\to M=G/_{\textstyle K}$ is 
realized by the diagram
$$
G=\a\times\KK\times\EK\stackrel{\pi}{\longrightarrow}\p=\a\times\EK
$$
where
$$
\pi(a,k,l)=(a,-\sinh(a)k+\cosh(a)l)
$$
with $\sinh(a) \stackrel{\mbox{def.}}{=}\frac{1}{2}\left( \exp\rho(a) - 
\exp\rho(-a)  \right)$ and $\cosh(a)\stackrel{\mbox{def.}}{=}\frac{1}{2}\left( 
\exp\rho(a) +  \exp \rho(-a)\right)$. The map $\gamma:\p \to G$ given 
by  
$$
\gamma(a,l)=(a,\sinh(a)l,\cosh(a)l)
$$
defines a global section of the principal bundle $(M,G,\pi)$. The action 
of $G$ on $M = \p$ reads as follows~:
$$
(\alpha,\kappa,\lambda)(a,l)=(a+\alpha,
\cosh(a+\alpha)\lambda-\sinh(a+\alpha)\kappa+l).
$$
\item[(ii)] The identification $M = \p$ defines a global Darboux 
chart $(M,\omega) \stackrel{\sim}{\to} (\p, \Omega)$. In this chart, 
the symmetries are given~by~: 
$$
s_{(a,l)}(a',l')=(2a-a',2\cosh(a-a')l-l').
$$
If $o=\pi(0,0,0)=(0,0)$, the midpoint map $M\to M:x\to \frac{x}{2} \quad 
(s_{\frac{x}{2}} o=x)$ is given (locally) by~:
$$
\frac{(a,l)}{2}=\frac{1}{2}\left( a, (\cosh\frac{1}{2}a)^{-1}l\right).
$$
\end{enumerate}
\end{prop}

The following map will be important while defining the deformed product.

\begin{dfn}\label{TWISTINGMAP} Let $(\g, \sigma, \Omega = \delta \xi)$ be 
an ESET (cf. Definition~\ref{ESET}). We denote by $\zeta$ the ``pairing"~:
$$
\zeta~: \a\times\EK\to \R
$$
defined by 
$$
\zeta(a,l)=\xi(\sinh(a)l).
$$
The formula 
$$
\omega(\phi(a),l) = \zeta(a,l)
$$
defines a ``linearization'' map $\phi : \a \to \a$, called the {\bf twisting 
map}. 
\end{dfn}

\begin{prop}\label{TWISTDIFFEO} Let $(M,\omega,s)$ be an elementary 
solvable symplectic symmetric space. Set $\omega = \Omega = \delta \xi$ 
($\xi\in \KK^*$) (see~\pref{SPACES} $(ii)$). Let $\a \times \EK 
\stackrel{\omega}{\to} \R$ be the symplectic pairing~:  
$$
\omega(a,l)=\Omega(a,l)=-\xi[a,l].
$$
Then, the midpoint map $M\to M:x\to \frac{x}{2}\quad(s_{\frac{x}{2}}o=x)$ 
is globally defined provided the twisting map $\phi : \a \to \a$ is a 
diffeomorphism. 
\end{prop}
\Pf
One has a linear isomorphism $\a\stackrel{\Omega}{\to}\EK^\star$ 
defined by $<\Omega(a),l>=\omega(a,l)$. Also, since $\zeta~: 
\a\times\EK\to \R$ is linear in its $\EK$-variable, one has a map 
$\a\stackrel{z}{\to}\EK^\star$ defined by $<z(a),l>=\zeta(a,l)$. One has 
 $<{z_\star}_a(A),l> = 
\xi(\rho(A)\cosh(a)l)=-\Omega(A,\cosh(a)l)$. \pref{TWISTDIFFEO} is now 
clear regarding the expression of the midpoint map given in \pref{SPACES}.
\EPf

\begin{rmks}\label{JAC}{\rm 
\begin{enumerate}
\item[(i)] This sufficient condition for existence of a global midpoint 
map is necessary as well, as it will be proven in Section~\ref{WKBQ} 
(Corollary~\ref{TMD}).
\item[(ii)] It is easy to determine the {\bf Jacobian}, 
$\mbox{Jac}_{\phi} (a)$, of the twisting diffeomorphism at point $a\in\a$. 
Indeed, identifying the tangent space $T_{a}(\a)$ with $\a$, one has 
for all $A\in\a$ and $l\in\EK$,
$$
\frac{d}{dt}|_{0}\Omega(\phi(a+tA),l)=\xi(\rho(A)\cosh(a)l)
$$
that is
$$
\Omega(\phi_{\star_{a}}(A),l)=\Omega(A,\cosh(a)l).
$$
Therefore, the linear map $\phi_{\star_{a}}~:\a\to\a$ is ``adjoint'' 
to the linear transformation $\cosh(a)|_{\EK}~:\EK\to\EK$. Hence
$$
\mbox{Jac}_{\phi}(a)=\det\left(\cosh(a)|_{\EK}\right).
$$
Note that $\cosh(a)|_{\EK}$ is indeed $\EK$-valued since $\rho( a)^{2} 
(\EK) \subset [\a,\KK] \subset \EK$. Observe also that in the case where 
$\g$ is nilpotent, one has $\mbox{Jac}_{\phi}(a)=1\,\forall a\in\a$.
\end{enumerate}
}\end{rmks}

\begin{dfn}\label{JC} 
Let $(\g=\b\oplus\a,\sigma, \Omega)$ be an ESET with splitting 
endomorphism $\rho : \a \to End(\b)$. For all $a\in\a$, write 
$$\rho(a) = \rho_N(a) + \rho_S(a)$$ for the Jordan-Chevalley decomposition 
of the (complex linearly extended) endomorphism $\rho(a):\b^c\to\b^c$, 
with $\b^c = \b \otimes \C$, and where $\rho_N(a)$ (respectively 
$\rho_S(a)$) denotes the nilpotent (respectively semisimple) part of 
the endomorphism $\rho(a)$.The ESET is said to be {\bf standard} if 
there exists $a^N$ and $a^S$ in $\a$ such that 
$$\begin{array}{ccl}
\rho_N(a) & = & \rho(a^N) \quad \mbox{and} \\
\rho_S(a) & = & \rho(a^S).
\end{array}$$
\end{dfn}

\begin{prop}\label{SESET}
Every elementary solvable exact triple is a sub-triple of a standard one. 
\end{prop}

\Pf
Let $(\g=\b\oplus\a,\sigma, \Omega)$ denote our starting exact triple. 
The map $\rho:\a\to End(\b)$ being injective (because $\Omega$ is 
nondegenerate), we may identify $\a$ with its image~: $\a=\rho(\a)$. Let 
$\Sigma:End(\b)\to End(\b)$ be the automorphism induced by the conjugaison 
with respect to the involution $\sigma|_{\b}\in GL(\b)$, i.e. $\Sigma = 
Ad( \sigma|_{\b})$. The automorphism $\Sigma$ is involutive and preserves 
the canonical Levi decomposition $End(\b)={\cal Z}\oplus sl(\b)$, where 
${\cal Z}$ denotes the center of $End(\b)$. Now, writing the element 
$a=\rho(a)\in\a$ as $a=a_Z+a_0$ with respect to this decomposition, one 
has~: $\Sigma(a)=a_Z+\Sigma(a_0)=-a=-a_Z-a_0$ because the endomorphisms 
$a$ and $\sigma|_{\b}$ anticommute. Hence $\Sigma(a_0)=-2a_Z-a_0$ and therefore 
$a_Z=0$. So, $\a$ actually lies in the semisimple part $sl(\b)$. For any 
$x\in sl(\b)$, we denote by $x=x^S+x^N,\quad x^S,x^N\in sl(\b)$ its 
abstract Jordan-Chevalley decomposition. Observe that, writing 
$sl(\b)=sl_+\oplus sl_-$ for the decomposition in $(\pm 1)$-$\Sigma$- 
eigenspaces, one has ~: $\a\subset sl_-$. Also, ${\cal N}:=\{ 
a^N\}_{a\in\a}$ is an Abelian subalgebra in $sl_-$ commuting with $\a$. 
One therefore may define the Abelian subalgebra in $sl_-$~:
$$
\a'=\a+{\cal N}.
$$
Canonically attached to $\a'$, one has the homomorphism~: $\rho':\a'\to 
End(\b)$ which anticommutes with $\sigma|_{\b}$, in particular ~: 
$\rho'(\a)\EK = \KK$. 

Let $\xi\in\KK^\star$ be the element whose coboundary defines the symplectic 
structure on $\p=\a\oplus \EK$ and denote by $\r\subset\a'\oplus \EK$ the 
radical of the coboundary $\delta\xi:\bigwedge^2(\a'\times_{\rho'}\b) 
\to \R$. Observe that $\r\subset\a'$. Hence, 
since $\a'\oplus\EK=\p\oplus\r$, one has $\a'=\a\oplus\r$. At last, 
let us denote by $\Lz(\g)$ the center of $\g$.\\

\noindent{\bf Case 1}~: $\Lz(\g)=0$.\\
Then $\g$ is a transvection algebra (cf.~\rref{DIM(Z)}). We form $\EK'=\EK \oplus\r^\star$, 
set $\b'=\KK\oplus \EK'$, and extend $\rho'$ to a homomorphism 
$\rho':\a'\to End(\b')$ as $\rho'(\a')\r^\star=0$. This homomorphism 
anticommutes with the involution $\sigma'|_{\b'}=id_\KK\oplus(-id_{\EK'})$, 
hence one has an involutive Lie algebra~: $(\a'\times_{\rho'} 
\b',(-id_{\a'}) \oplus \sigma'|_{\b'})$. Now, we define a symplectic structure 
$\Omega'$ on $\p'=\a'\oplus \EK'$~by~:
$$
\begin{array}{c}
\Omega'|_{\p\times\p}=\Omega|_{\p\times\p}\\
\Omega'|_{\p\times(\r\oplus\r^\star)}=0\\
\Omega'(r^*,r)=<r^*,r> \quad \forall r \in \r ; r^* \in \r^\star.
\end{array}
$$
The symplectic structure $\Omega'$ turns out to be $\KK$-invariant. 
Indeed, one has~:
$
\left[\KK,\r^\star\right]=\left[\left[\a,\EK\right],\r^\star\right]=0 
\mbox{ (by Jacobi) and }\left[\KK,\r\right]\subset\EK.
$
Now, considering a one-dimensional central extension of 
$\a'\times_{\rho'}\b'$, one gets the desired standard exact triple.\\

\noindent{\bf Case 2}~: $\Lz(\g)\neq0$.\\
Define $E\in\Lz(\g)\subset\KK$ by $\xi(E)=1$. Set $\EK'=\EK 
\oplus\r^\star$, $\b'=\KK\oplus \EK'$, and 
endow the vector space $\g'=\b'\oplus \a'$ with the skewsymmetric 
bracket defined by~:
$$
\begin{array}{c}
\left[\, , \, \right]_{{\g'}}|_{\g\times\g}=\left[\, , \, 
\right]_{\g}\\
\left[r^{*},r\right]=<r^{*},r>E \quad \forall r \in \r ; r^* \in 
\r^\star\\
\left[r,b\right]=\rho'(r)b \quad \forall r \in \r ; 
b\in\b\subset\b',
\end{array}
$$
the other brackets being zero. Endowed with this structure $\g'$ 
turns out to be a Lie algebra. Moreover the coboundary 
$\Omega'=\delta\xi:\bigwedge^{2}(\g')\to\R$ restricts to 
$\p'=\a'\oplus\EK'$ as a symplectic 
structure. Therefore the triple 
$(\g', id_{\KK}\oplus(-id_{\p'}),\Omega')$ is the desired 
standard exact triple.
\EPf

\begin{rmk}
{\rm 
Observe that the exact triples constructed in cases 1 and 2 above are both such 
that $$\dim\z(\g)=1.$$ This ensures they are one-dimensional central extensions of 
(transvection) symplectic triples.
}
\end{rmk}

Led by the example in the preceding section, one observes
\begin{prop}\label{PHASESET}
Let $(M=\a\times\EK,\omega,s)$ be an elementary solvable symplectic 
symmetric space. Then, the function 
$$
\tilde{u}:M\times M\to\KK
$$
defined by
$$
\tilde{u}((a,l),(a',l'))=\sinh(a')l-\sinh(a)l'
$$
satisfies (locally)
$$
\tilde{u}(x,y)=-\tilde{u}(y,x)=-\tilde{u}(x,s_x(y))=-\tilde{u}(x,s_{\frac{x}{2}}(y))\qquad\forall 
x,y\in M
$$
and is $\exp(\KK)$-invariant. In particular, for every choice of an 
element $\xi\in\KK^\star$, the function 
$$
u=\xi\circ\tilde{u}
$$
is an admissible two-point function on $M$ (Definition~\ref{TPF}) provided the midpoint map 
is globally defined. The associated admissible three-point function 
(Proposition~\ref{MIDPOINT}) has 
the form~:  
$$
S((a_1,l_1),(a_2,l_2),(a_3,l_3))=\xi\left(\oint_{1,2,3}\sinh(a_1-a_2)l_3\right)
$$
where $\oint_{1,2,3}$ stands for the cyclic summation.
\end{prop}

One sees here why it is important to deal with an exact triple rather 
than the transvection triple. Indeed, one would need, while proceeding 
a stationary phase method for the deformed product (cf.~Dirac's condition 
in definition~\ref{WKB}), to obtain the Poisson bracket associated to the 
the symplectic structure $\omega$ as the first order expansion term. The 
critical point analysis tells us that the above expression of $S$ yields 
a first order term equal to $\delta\xi$.\\  

\section{Obtaining the oscillating kernel via star representation theoretical 
methods}\label{REPTH}

In this section, we present a heuristic way to derive oscillatory 
integral formulae for our strict deformation quantizations (see 
\tref{WKBTHM}). For the sake of simplicity, we will only treat the case of 
the two-dimensional example $(\R^2, \omega^0,s)$, already mentioned in 
Example~\ref{EX2DIM}. The general case will be treated later. The 
proof of \tref{WKBTHM} does not depend on the present section, but we believe 
that it clarifies the ``mysterious'' formula (\ref{POMME}) of the deformed 
product. This section, together with Section~\ref{WEYLR}, suggests that our 
construction should generalize to a much wider situation than solvable 
symmetric space. 

Following a classical result of Kostant, the symplectic homogeneous 
space $(M=G/_{\textstyle K},\omega)$ is an equivariant symplectic covering 
of a coadjoint orbit in $\g^\star$. From Proposition~\ref{SPACES}, one gets 
a global section $\gamma : \p =  \a \times \EK \to G=\a \times \KK \times 
\EK$. This yields a global Darboux chart on the coadjoint orbit ${\cal 
O} = Ad^{\star}(G)\xi\subset\g^{\star}$~:   
$$
\p\to{\cal O}: x\to Ad^{\star}(\gamma(x)) \xi.
$$
In particular, the Hamiltonian function $\lambda_{X}\in C^{\infty}(\p)$ 
associated to the infinitesimal action of $X \in \g$  (see \sref{SSS}) is 
given by~: 
\begin{eqnarray*}
\lambda_{X}(x) & = & <Ad(\gamma(x))\xi,X>\\
 & = & <\xi,Ad(\gamma(x)^{-1})X>\\
 & = & <\xi, \rho(X_{\a})l-R(-a)X_{\b}>
\end{eqnarray*}
where $x=(a,l)$ and $X=X_{\a}+X_{\b}\in\g=\a\times\b$. From this last 
expression, one obtains for $X\in\p$~: 

\begin{eqnarray}\label{DERIVEES}
\partial_{l}^{\beta}\lambda_{X} & = & 0 \quad (\beta\geq 2)\\
\partial_{a}^{\alpha}\partial_{l}^{\beta}\lambda_{X} & = & 0 \quad 
(\alpha,\beta\geq 1) \nonumber \\
\partial_{l}\lambda_{X}(x) & = & <\xi, \rho(X_{\a})\partial_{l}> 
\nonumber \\
\partial_{a}^{\alpha}\lambda_{X}(x) & = & \cosh(a) <\xi, 
\rho(\partial_{a})^{\alpha} X> \quad (\alpha \mbox{ odd }, \alpha\geq 
1), \nonumber
\end{eqnarray}
where $\partial_a$ (respectively $\partial_l$) is an element of $\a$ 
(respectively $\EK$) thought of as a constant vector field on $\p$.
In particular, if $\star_{\nu}^{M}$ denotes the standard Moyal star 
product on $(\p,\Omega)$ (cf.~formula~(\ref{MOYAL}) in Section~\ref{WEYLR}), 
one gets 
$$
[\lambda_{X}, \lambda_{Y}]_{\sm} \stackrel{\mbox{def.}}{=} \lambda_{X} \sm 
\lambda_{Y} -\lambda_{Y}\sm\lambda_{X}=2\nu\{\lambda_{X}, 
\lambda_{Y}\}\quad \forall X,Y\in\g. 
$$
One refers to this last property as the {\bf $\g$-covariance} of the 
Moyal star product (\cite{AC2}). The covariance property allows us 
to define a representation of $\g$ on the space $C^{\infty}({\cal 
O})[[\nu]]$ of formal power series in the parameter $\nu$ with 
coefficients in $C^{\infty}({\cal O})$~:

\begin{eqnarray*}
\rho_{\nu}: \g\to End(C^{\infty}({\cal O})[[\nu]])~: \\
\rho_{\nu}(X) u =\frac{1}{2\nu}[\lambda_{X}, u]_{\sm}.
\end{eqnarray*}
Formula~(\ref{MOYAL}) in Section~\ref{WEYLR} yields, using 
(\ref{DERIVEES}),  

\begin{eqnarray*}
(\rho_{\nu}(X)u)(x) & = & \{\lambda_{X}, 
u\}-\frac{1}{\nu}\sum_{k=1}^{\infty}
\frac{\nu^{2k+1}}{(2k+1)!}
\cosh(a)<\xi,\rho(\partial_{a})^{2k+1}X>\partial_{l}^{2k+1}u \\
 & = & <\xi, \rho(X_{\a})\partial_{l}> \partial_{a}u 
 -\frac{\cosh(a)}{\nu}\sum_{k=0}^{\infty}\frac{\nu^{2k+1}}{(2k+1)!}
 <\xi,\rho(\partial_{a})^{2k+1}X>\partial_{l}^{2k+1}u .           
\end{eqnarray*} 
A partial Fourier transform in the $\EK$-variables allows to 
interpret $\rho_{\nu}$ as a ``multiplicative representation''. 
Indeed, setting 
$$
F(u)(a,\alpha)=\hat{u}(a,\alpha)=\int_{\EK}e^{-i\Omega(\alpha,l)}u(a,l)\,dl,   
$$
one gets~:
$$
{\hat{\rho}}_{\nu}(X).\hat{u}\stackrel{\mbox{def.}}{=}F(\rho_{\nu}(X)u)= 
X_{\a}.\hat{u}-\frac{\cosh(a)}{\nu}\sum_{k=0}^{\infty}\frac{(i\alpha
\nu)^{2k+1}}{(2k+1)!}<\xi,\rho(\partial_{a})^{2k+1}X>\hat{u}         
$$
for all $X=X_{\a}+X_{\EK}\in\p$. That is, setting 
$\nu=\frac{\hbar}{2i}$, for all $X\in\p$~:
$$
\hat{\rho}_{\frac{\hbar}{2i}}(X).\hat{u}=X_{\a}.\hat{u}-\frac{2i\cosh(a)}{\hbar}
<\xi,\sinh(\frac{\hbar\alpha}{2})X>\hat{u}
$$
or 
$$
\hat{\rho}_{\frac{\hbar}{2i}}(X).\hat{u}=X_{\a}.\hat{u}+c_{\hbar}(X)\hat{u}
$$
with $c_{\hbar}~:\p\to C^{\infty}(\R^{2})$ defined by 
$$
c_{\hbar}(X)(a,\alpha) \stackrel{\mbox{def.}}{=} -\frac{2i\cosh(a)}{\hbar}<\xi,
\sinh(\frac{\hbar\alpha}{2})X>.
$$
The expression of the cocycle $c_{\hbar}$ is very similar to the one 
of the ``twisting map'' $\zeta$ associated to the three-point 
function $S$ (cf.~\pref{TWISTDIFFEO}). From this observation, we now derive a 
commutative multiplicative law on the functions on $\R^{2}=\{ (a,\alpha) \}$ 
for which representation $\hat{\rho}$ acts by derivations. We start by 
observing that, for all $\partial_{l}\in\EK$, one has 
$\partial_{l}\zeta(\alpha,l)=<\xi,\sinh(\alpha)\partial_{l}>$.
Hence, setting 
$$
\zeta_{\hbar}(\alpha,l) \stackrel{\mbox{def.}}{=} \frac{2}{\hbar} \zeta( 
\frac{\hbar}{2} \alpha,l)    
$$ 
and 
$$\mu=-\frac{1}{i\cosh(a)},$$ 
one gets $\partial_{l}\zeta_{\hbar} 
(\alpha,l) = \mu\, c_{\hbar}(\partial_{l})(a,\alpha)$. Therefore, 
defining formally 
$$
(\Lz_{\hbar}u)(a,\alpha) \stackrel{\mbox{def.}}{=} \int_{\EK}e^{-i 
\zeta_{\hbar}(\alpha,l)} u(a,l)\,dl 
$$
for (reasonable) $u\in C^{\infty}(\p)$, an integration by parts argument 
leads us to $\Lz_{\hbar}(\partial_{l}u)=i\mu\, c_{\hbar}(\partial_{l}) 
\Lz_{\hbar} (u)$ or  
\begin{equation}\label{DELEL}
\Lz_{\hbar}^{-1}(c_{\hbar}(\partial_{l}) 
f)=\frac{1}{i\mu}\partial_{l}\Lz_{\hbar}^{-1}(f)
\end{equation}
whenever this last formula makes sense for $f\in C^{\infty}(\R^{2})$.
When well-defined, the following (commutative) product
$$
f \bullet_{\hbar} g \stackrel{\mbox{def.}}{=} 
\Lz_{\hbar}(\Lz_{\hbar}^{-1}f.\Lz_{\hbar}^{-1}g)
$$ 
is then invariant under the representation $\hat{\rho}$. Indeed, 
since the integral in the definition of $\Lz_{\hbar}$ is only taken 
on the $\EK$-variable, the product $\bullet_{\hbar}$ is invariant under 
$\a$ (observe that $c_{\hbar}(\a)=0$). The algebra $\g$ being generated 
by $\p$, it is therefore sufficient to prove $c_{\hbar}(\partial_{l}) f 
\bullet_{\hbar} g =(c_{\hbar}(\partial_{l}) f)\bullet_{\hbar} g + f 
\bullet_{\hbar} (c_{\hbar}(\partial_{l}) g)$. The RHS is actually 
$\frac{1}{i\mu}\Lz_{\hbar} ((\partial_{l}\Lz_{\hbar}^{-1}f).\Lz_{\hbar}^{-1}
g+\Lz_{\hbar}^{-1}f.\partial_{l}\Lz_{\hbar}^{-1}g)$
using (\ref{DELEL}) and the fact that multiplication by $\mu$ commutes 
with the transformation $\Lz_{\hbar}$. One therefore has 
$\frac{1}{i\mu}\Lz_{\hbar}\partial_{l}(\Lz_{\hbar}^{-1}f.\Lz_{\hbar}^{-1}g)$ or 
$\frac{1}{i\mu}\Lz_{\hbar}\partial_{l}\Lz_{\hbar}^{-1}(f \bullet_{\hbar} 
g )$ which equals $c_{\hbar}(\partial_{l}) f \bullet_{\hbar} g $ 
using (\ref{DELEL}) again.\\
If one interprets the commutative product $\bullet_{\hbar}$ as the 
underlying product to the algebra of functions on a commutative 
$\hbar$-dependent manifold, say $M_{\hbar}$, its invariance under 
$\hat{\rho}$ tells us that $\g$ is realized via $\hat{\rho}$  as a 
subalgebra of tangent vector fields over $M_{\hbar}$. It turns out 
that this action of $\g$ on $M_{\hbar}$ is equivalent to the one of $\g$ 
on $M$. Indeed for all $X\in \g$, one has~: 
\begin{equation}\label{MHBAR}
\Lz_{\hbar}^{-1}\circ \hat{\rho}_{\frac{\hbar}{2i}}(X) \circ\Lz_{\hbar}=X^{\star}.
\end{equation}
Again, it is sufficient to prove this last formula only for 
$X=\partial_{l}\in\EK$. In this case, formula (\ref{DELEL}) identifies 
the LHS action on $u\in C^{\infty}(\p)$ to $\frac{1}{i\mu}\partial_{l}u$ 
which is precisely $\partial_{l}^{\star}u$.\\
Formula (\ref{MHBAR}) leads us to consider the transformation $T_{\hbar}$ 
defined by 
$$
T_{\hbar}=F^{-1}\circ\Lz_{\hbar}
$$ 
which intertwines our $\star$-representation $\hat{\rho}_{\frac{\hbar}{2i}}$ 
with the infinitesimal (coadjoint) action of the transvection algebra. 
In other words, the formal product on $C^{\infty}(\p)[[\hbar]]$ defined by 
$$
u\star_{\hbar} v \stackrel{\mbox{def.}}{=} T_{\hbar}^{-1} (T_{\hbar} u 
\star_{\frac{\hbar}{2i}}^{M} T_{\hbar}v)    
$$
is invariant under the coadjoint action of $G$ on ${\cal O}=\p$.

\section{WKB-Quantization}\label{WKBQ}

\begin{prop}\label{GROWTH-P}
Let $(M,\omega,s)$ be an elementary solvable symplectic symmetric 
space. Let $\phi : \a \to \a$ be the associated twisting map 
(cf.~\dref{TWISTINGMAP}). Assume that $(M,\omega,s)$ admits a globally 
defined midpoint map (cf.~\pref{SPACES} $(ii)$). Let $\|,\|$ be a norm 
on the vector space $\a$. Then 
$$
\liminf_{a \to \infty} \frac{\| \phi (a) \|}{\|a\|} > 0.
$$
\end{prop}
Before passing to the proof, we observe

\begin{lem}\label{GROWTH}
Let $V$ be a finite dimensional real vector space. Let $\a$ be a non 
trivial Abelian subalgebra of $End(V)$
such that 
\begin{enumerate}
\item[(i)] if the Jordan-Chevalley decomposition (cf.~\dref{JC}) of 
$a\in\a$ writes~: $a=a^S+a^N$ then $$a^S, a^N\in \a;$$ 
\item[(ii)] every non zero (complex) eigenvalue of $a^S$ has a non trivial 
real part.
\end{enumerate}
Consider the function $\a\to End(V): A \to\sinh(a)$. Then, 
$$
\liminf_{a\to\infty}\frac{||\sinh(a)||_{op}}{||a||_{op}}>0
$$
where $||\, , \, ||_{op}$ denotes the operator norm on $End(V)$ defined 
with respect to any norm compatible with the vector space structure on $V$.
\end{lem}

\Pf
Let us first extend the action of $\a$ on $V$ $\C$-linearly to a 
(complex) action of $\a$ on $V^c=V\otimes\C$. Denote by 
$V^c=\oplus_{\lambda\in\Phi}V_\lambda$ the weight space decomposition 
with respect to the Abelian subalgebra of semisimple elements~: $\s=\{ 
a^S \}_{a\in\a}$, where $\Phi$ denotes the set of weights of the action 
of $\s$. Observe that, setting ${\cal N}=\{a^N\}_{a\in\a}$, one has ${\cal 
N}.V_\lambda\subset V_\lambda$ for all $\lambda$ in $\Phi$. In each 
$V_\lambda$, fix a basis $e_\lambda$ for which the matrix 
associated to any $n\in{\cal N}$ is upper triangular. In particular, in 
the basis $e=\{ e_\lambda \}_{\lambda\in \Phi}$, the matrix associated to 
any $a\in\a$ is upper triangular with elements $\lambda(a^S)$ on the 
principal diagonal. Let $\tilde{\a}$ be a maximal Abelian subalgebra of 
$End(V^c)$ which contains $\a$. Maximality implies 
that, similarly to $\a$, one has a decomposition in semisimple and nilpotent parts~:
$\tilde{\a}=\tilde{\s}\oplus\tilde{{\cal N}}$. Note that 
$\sinh(\a)\subset\sinh(	\tilde{\a})\subset \tilde{\a}$.\\
On $ \tilde{\a}$, we consider the norm~:
$$
||\tilde{a}=\tilde{s}+\tilde{n}||\stackrel{\mbox{def.}}{=}\mbox{max}\{||
\tilde{s}||_{op}, ||\tilde{n}||_{op}\}\quad 
(\tilde{s}\in\tilde{\s},\tilde{n}\in\tilde{\cal N}),
$$
where $||\,||_{op}$ denotes the operator norm on $End(V^c)$ 
associated to the choice of a norm on $V^c$. Now, we observe that for all 
sequences $\{s_k\}_{k\in\N}\subset \s$, one has 
$$\lim_{k\to \infty} \frac{|| \sinh(s_k)||}{||s_k||}\neq 0.
$$ 
The proof of the last assertion is divided into the three following steps. 

First, for any fixed $k\in \N$, one can find a weight 
$\lambda^{(k)}\in\Phi$ which realizes $|\lambda^{(k)}(s_k)|=||s_k||_{op}$. 
In particular, 
$$
\frac{||\sinh(s_k) ||}{||s_k||} \geq 
\frac{|\sinh(\lambda^{(k)}(s_k))|}{|\lambda^{(k)}(s_k)|}.  
$$ 
Assume by the absurd that 
$$\lim_{k\to\infty}\frac{||\sinh(s_k)||}{||s_k||}=0.$$
Since $\Phi$ is a finite set, one can then find a weight $\lambda\in\Phi$ and a 
partial sequence $\{s'_l\}\subset\{s_k\}$ such that 

\begin{equation}\label{inter1}
\lim_{l \to \infty}\frac{|\sinh (\lambda (s'_l))|}{|\lambda (s'_l)|}=0. 
\end{equation}

Secondly, observe that every weight $\lambda\in\Phi$ writes 
$\lambda=\alpha_\lambda+i\beta_\lambda$ with $\alpha_\lambda$ and 
$\beta_\lambda$ in the real vector space $\s^\star$ with the property 
that if $s\in \ker\alpha_\lambda$ then $\lambda(s)=\beta_\lambda(s)=0$ 
(hypothesis (ii)). Hence $\ker \alpha_{\lambda}\subset \ker 
\beta_{\lambda}$, and there exists a real $t_{\lambda}\in\R$ 
such that $\lambda=(1+it_{\lambda})\alpha_{\lambda}$. 

Thirdly, for all $z=x+iy\in\C$, one has $|\sinh(z)|^{2} = \sinh (x)^2 + 
\sin(y)^{2}$. Hence writing $\lambda(s'_{l})=x_l+iy_{l}$, the 
relation (\ref{inter1}) yields 

\begin{equation}\label{inter2}
\lim_{l\to\infty}\frac{\sinh(x_{l}^{2})}{x_{l}^{2}+y_{l}^{2}}=0.
\end{equation}
But, by our preceding observation, there exists a real number $t$ such 
that $y_{l}=tx_{l}$. Hence $|\lambda (s'_{\lambda})|^{2} = 
(1 + t^{2})x_{l}^{2}$. This, combined with (\ref{inter2}), yields a real 
sequence $x_{l}$ such that $\lim_{l \to \infty} \frac{ 
\sinh(x_{l}^{2})}{x_{l}^{2}} = 0$, and we reach a contradiction.

Suppose, that $\{a_k\}_{k\in\N}$ is a sequence of elements in $\a$ for 
which the sequence $||a_k||$ diverges to infinity and such that 
$\lim_{k\to\infty}\frac{||\sinh(a_k)||_{op}}{||a_k||_{op}}=0$.
Every two norms on $ \a$ being equivalent, one has~: 
$\lim_{k\to\infty}\frac{||\sinh(a_k)||}{||a_k||}=0$. Consider the 
following sequences partitioning 
$\{a_k\}$~:
$$
\begin{array}{ccc}
\{\nu_p\}& = & \{a_k=a_k^S+a_k^N\mbox{ such that 
}||a_k^N||>||a_k^S||\}\\
\{\sigma_q\}& = & \{a_k=a_k^S+a_k^N\mbox{ such that 
}||a_k^N||\leq||a_k^S||\}.
\end{array}
$$
One has 
$\lim_{q\to\infty}\frac{||\sinh(\sigma_q)||}{||\sigma_q||}=0$. And, since 
$$
\frac{||\sinh(\sigma_q)||}{||\sigma_q||}
=\frac{\mbox{max}\{||(\sinh(\sigma_q))^S||,
||(\sinh(\sigma_q))^N||\}}{||\sigma^S_q||}\geq 
\frac{||(\sinh(\sigma_q))^S||}{||\sigma^S_q||}
=\frac{||\sinh(\sigma^S_q)||}{||
\sigma^S_q||},
$$ 
our previous observation implies that $\{\sigma_{q}\}=\emptyset$, 
i.e.~that $\{a_k\}=\{\nu_p\}$. As above, one gets $\lim_k \frac{||( 
\sinh(a_k))^N||}{|| a^N_k||}=0$, that is, 
\begin{equation}\label{1}
\lim_{k\to\infty}\frac{1}{||a_k^N||}\left\{\cosh(a^S_k)\sinh(a^N_k)
+\sinh(a^S_k)(\cosh(a^N_k)-{\bf 1})\right\}=0 .
\end{equation}

For all $a\in\a$, the endomorphisms $\sinh(a^N)$ and $\cosh(a^N)-{\bf 1}$ 
are linearly independent in $End(V^c)$. Indeed, those are sums of powers 
of the upper triangular nilpotent matrix $a^N$. Hence, none of the terms 
of $\cosh(a^N)-{\bf 1}$ can cancel the term $a^N$ in $\sinh(a^N)$. Therefore, 
since $\cosh(a^S)$ as well as $\sinh(a^S)$ act diagonally, the terms occurring 
in the bracket of expression~(\ref{1}) are linearly independent. In 
particular, $\lim_k\cosh(a^S_k)\frac{1}{||a_k^N||} \sinh(a^N_k)=0$. Therefore, 
$\cosh(a^S_k)$ being invertible for all $k$, one gets 
$$\lim_k\frac{1}{||a_k^N||} \sinh(a^N_k)=0.$$ Again, since the matrix $a_k^N$ is upper 
triangular, one gets $\lim_k \frac{1}{||a_k^N||}a_k^N=0$, a contradiction.
\EPf

{\em Proof of \pref{GROWTH-P}}. Proposition~\ref{SESET} 
reduces to prove the assertion only for standard ESET's. Indeed, the 
twisting map associated to any ESET can be realized as the restriction 
of the twisting map of a standard ESET to a linear subspace. We therefore 
assume our ESET $(\g,\sigma,\Omega)$ to be standard.\\    

Let us define $C_{sl_-}(\a)$ to be the intersection of the centralizer 
algebra of $\a=\rho(\a)$ in $End(\b)$ and $sl_-$ (cf.~Proof of \pref{SESET}). 
For $X\in C_{sl_-}(\a)$, we set 
$$
||X||\stackrel{\mbox{def.}}{=}\sup_{||l||\leq1}\{|\xi(X.l)|\}
$$ where $l\in\EK=\b\cap\p$ and $\delta\xi=\Omega\quad(\xi\in\k^\star)$.
This defines a norm on $C_{sl_-}(\a)$. Indeed, if $X\in C_{sl_-}(\a)$ 
is such that $\xi(X.\EK)=0$, one has $\xi(X. \left[\a,\k\right]) = 
\xi \left[\a,X.\k\right] =\Omega(\a,X.\k)=0$, hence $X.\k=0$.
Also, $0=X.\k=X.\left[\a,\EK\right]=\left[\a,X.\EK\right]$, hence 
$X.\EK\subset\z(\g)$. By item (ii) of \rref{DIM(Z)}, one has either 
$X.\EK=0$ or $X.\EK\subset\R.E$ with $\xi(E)=1$, which implies $X.\EK=0$ 
too since one assumed $\xi(X.\EK)=0$. Now, observe that $\sinh(\a)\subset 
C_{sl_-}(\a)$, and, that for all $a\in\a$, one has~:
$$
||\sinh(a)||=\sup_{||l||\leq1}\{|\xi(\sinh(a)l)|\}=\sup_{||l||\leq1}\{|
\Omega(\phi(a),l)|\}=||\phi(a)||.
$$
Therefore, the norms $||\, , \, ||_{op}$ and $||\, , \, ||$ being equivalent, 
\lref{GROWTH} implies \pref{GROWTH-P}.
\EPf

We now have the following geometric property.

\begin{cor}\label{TMD}Let $(M,\omega,s)$ be an elementary solvable symplectic 
symmetric space. Then $(M,\omega,s)$ admits a globally defined 
midpoint map $x \mapsto \frac{x}{2}$ (if and) only if the twisting 
map $\phi : \a \to \a$ (cf.~\dref{TWISTINGMAP}) is a global 
diffeomorphism.
\end{cor}

\Pf With the same notations as in the proof of \pref{TWISTDIFFEO}, one 
has a map $\phi : \a \to \a$ defined by $\phi = \Omega \circ z$. The 
argument used in the same proof tells us that $\phi$ is a local 
diffeomorphism. Now, \pref{GROWTH-P} implies that the map $\phi : \a \to 
\a$ is  proper. Indeed, $\phi$ being continuous, one just needs to show that 
the inverse image of a ball is bounded. So, let $\{x_k\}\subset\a$ be such 
that
$$
||\phi(x_k)||<R.
$$
Since $0<c<\liminf\frac{||\phi(x_k)||}{||x_k||}\leq\frac{R}{||x_k||}$, one 
gets $||x_k||\leq\frac{R}{c}$ as soon as $k$ is large enough.\\
Therefore, nonetheless the map $\phi$ is open it is also closed, hence 
surjective.
 The twisting map is therefore a covering 
map. Since the fundamental group of $\a$ is trivial, it is a 
diffeomorphism.
\EPf

\begin{ntt}
{\rm
Let $V$ be a finite dimensional real vector space. We denote by ${\cal 
S}(V)$ the space of Schwartz (complex-valued) functions on $V$. 
Accordingly, ${\cal S}'(V)$ stands for the space of tempered 
distributions on $V$. 
}
\end{ntt}

\begin{prop}\label{STAB}
Let $(M,\omega,s)$ be an elementary solvable symplectic symmetric 
space, with associated ESET $(\g = \a\oplus \b ,\sigma, \Omega)$ 
(cf.~\dref{ESET}). Assume that $(M,\omega,s)$ admits a globally defined 
midpoint map. Let $\phi :\a \to \a$ be the associated twisting 
diffeomorphism (cf.~\cref{TMD}). Then, one has

\begin{enumerate}
\item[(i)] $\phi^{\star}\s(\a)\subset\s(\a)$ and
\item[(ii)] $(\phi^{-1})^{\star}\s(\a)\subset\s'(\a)$.
\end{enumerate}
\end{prop}
We use the following lemmas.

\begin{lem}\label{COSECH}
Within the hypotheses and notations of \lref{GROWTH}, the function
$$\a\to End(V):a\to(\cosh(a))^{-1}$$ is well-defined as a tempered 
analytic function from $\a$ to the (real Banach) space $End(V)$.
The same holds for $a\to\tanh(a)$.
\end{lem}

\Pf
One has $\cosh(a)=\cosh(a^S+a^N)=\cosh(a^S)\left(\cosh(a^N)+\tanh(a^S)
\sinh(a^N)\right)$. In the basis $e$ considered in the proof of 
\lref{GROWTH}, $\cosh(a^S)$ and $\tanh(a^S)$ are diagonal matrices whose 
elements are in $\{\cosh(\lambda(a^S)), \tanh( \lambda(a^S))\}_{\lambda 
\in \Phi}$. Since each weight is of the form $\lambda=(1+it)\alpha$ 
(cf.~Proof of \lref{GROWTH}), $\cosh(a^S)^{-1}$ is a Schwartz function, 
and $\tanh(a^S)$ together with all its derivatives are bounded.

The factor $\cosh(a^N)+\tanh(a^S)\sinh(a^N)$ is an upper triangular matrix 
whose diagonal elements are equal to $1$. Moreover, it is a polynomial 
in the variable $a^N$. Therefore its inverse is a tempered function in 
the variable $a$. A similar argument (but simpler) yields the second part 
of the assertion. 
\EPf

\begin{lem}\label{SINH} With the hypotheses of \pref{GROWTH-P}, one has~: 
$$
\sinh(\phi^{-1}(a))=\rho(a)
$$
for all $a\in\a$.
\end{lem}
\Pf
One has, for all $k\in\k$ and $a,A\in\a$~: 

$$
\begin{array}{c}
\Omega(\sinh(\phi^{-1}(a))k,A)=\xi\left[\sinh(\phi^{-1}(a))k,A\right]=\xi
(\sinh(\phi^{-1}(a))\left[k,A\right])=\\
\zeta(\phi^{-1}(a),\left[k,A\right])=\Omega(\phi(\phi^{-1}(a)),\left[k,
A\right])=\Omega(a,\left[k,A\right])=\Omega(\left[a,k\right],A).
\end{array}
$$
Hence $\sinh(\phi^{-1}(a))k=\left[a,k\right]$. Moreover, 
$\sinh(\phi^{-1}(a))\left[k,A\right]=\left[\sinh(\phi^{-1}(a))k,A\right]=\left[
\left[a,k\right],A\right]=\left[a,\left[k,A\right]\right]$.
\EPf

\begin{lem}\label{BANACH}
Let $\b(E)$ be the space of bounded linear operators on the Banach 
space $E$. Let $a\in \b(E)$ be invertible and such that 
$||a^{-1}||\leq1$. Then, $||a^{2}||\geq||a||$.
\end{lem}
\Pf
One has $\sup_{x\in 
E,||x||\leq1}\{||a^{2}a^{-1}x||\}\leq\sup_{y\in 
E,||y||\leq||a^{-1}||}\{||a^{2}y||\}\leq\sup_{x\in 
E,||x||\leq1}\{||a^{2}x||\}$.
\EPf

\vspace{2mm}

{\em Proof of \pref{STAB}.} As in \pref{GROWTH-P}, we can assume our ESET 
to be standard. A function $u\in C^{\infty}(\a)$ is Schwartz if and only 
if, for all multiindex $\alpha$, and all positive integer $N$, one has 
$\sup_{a\in\a}\{||a||^{N}|(D^{\alpha}u)(a)|\}<\infty$. Consider 
$u\in\s(\a)$. Then, one has 
$$\sup_{a\in\a}\{||a||^{N}|(\phi^{\star} u) 
(a)|\} =\sup_{a\in\a}\{||\phi^{-1} (a)||^{N}|u(a)|\}.$$ 
By \pref{GROWTH-P}, one can find $r>0$ such that
$$
\sup_{a\in\a}\{ ||\phi^{-1} (a) ||^N|u(a)|\} \leq r 
\sup_{a\in\a}\{||a||^{N}|u(a)|\}. 
$$
Now, consider $A\in\a$ and let $$D_{A}(\phi^{\star}u)(a) =  
\frac{d}{dt}|_{0} \phi^{\star}u(a+tA).$$ In order to bound 
$\sup_{a\in\a}\{ ||a||^{N}|(D_{A}(\phi^{\star}u))(a)|\}$, the preceding 
argument leads us to look at $$<\phi_{\star\phi^{-1}(a)}(A),du|_{a}>,$$ that 
is, at $||\phi_{\star\phi^{-1}(a)}(A)||$. Since 
$$
\frac{d}{dt}|_{0} \Omega(\phi( \phi^{-1}(a)+tA),l) =\xi (\cosh ( 
\phi^{-1} (a)) \rho(A)l)
$$ (using the definition of $\phi$), one only needs to analyze the asymptotic 
behavior of $\cosh(\phi^{-1}(a))$. For this, observe that for some $m>0$ 
and for all $a\in\a$, one has~:   
$$
||a||^{-m}||(\cosh(\phi^{-1}(a)))^{-1}||= ||\phi(a')||^{-m}||
(\cosh(a')))^{-1}||\leq||a'||^{-m}||(\cosh(a')))^{-1}||,
$$ 
as soon as $||a'=\phi^{-1}(a)||$ is large enough. Therefore, using 
\lref{COSECH}, there exists $m>0$ such that
$$
\left|\left| 
||a||^{-m}(\cosh(\phi^{-1}(a)))^{-1} \right|\right|\leq1
$$ 
for $||a||$ large enough. Hence, \lref{BANACH} implies~: 

$$
||a||^{2m}||(\cosh(\phi^{-1}(a)))^{2}||
\geq||a||^{m}||\cosh(\phi^{-1}(a))||
$$ 
that is, by \lref{SINH}~:

$$
||\cosh(\phi^{-1}(a))||\leq||a||^{m}||{\bf 
1}+(\rho(a))^{2}||
$$ 
which has a polynomial growth. 

Now, we indicate how to prove (i) by induction over the order of 
derivation. Let $\nabla$ be the flat Euclidean connection on the vector 
space $\a$. For $u\in C^\infty(\a)$, let $\nabla^{(r)}u$ be the
symmetrization of the tensor field $\nabla^ru$. For the two first orders 
of derivation, one then has
$$
D_A(\phi^\star u)=<\nabla 
u\circ\phi,(\phi_\star(A))\circ\phi>=\phi^\star<\nabla u,\phi_\star A>
$$
and 
$$
D_BD_A(\phi^\star u)=\phi^\star(<<\nabla^{(2)} u,\phi_\star B>,\phi_\star A>
+<\nabla u, <\nabla(\phi_\star A),\phi_\star B>>).
$$
The functions (sections) $\nabla u$ and $\nabla^{(2)}u$ are Schwartz and 
we have seen that $\phi_\star A$ has polynomial growth.
So, in order to bound the second derivative $D_BD_A(\phi^\star u)$, one 
needs to control the asymptotic behavior of $\nabla(\phi_\star A)$, that is 
$\Omega(<\nabla(\phi_\star A),B>,l)$ for all $A,B\in\a, l\in\EK$. One has
\begin{eqnarray*}
\Omega(<\nabla(\phi_\star 
A),B>|_a,l)=\frac{d}{dt}|_0\Omega(\phi_{\star_{\phi^{-1}(a+tB)}}(A),l)=\\
\xi(\sinh(\phi^{-1}(a))(\cosh(\phi^{-1}(a)))^{-1}\rho(A)\rho(B)l)=\\
\xi(\rho(a)(\cosh(\phi^{-1}(a)))^{-1}\rho(A)\rho(B)l).
\end{eqnarray*}
We have seen previously that $(\cosh(\phi^{-1}(a)))^{-1}$ has polynomial
growth. Therefore,  $\nabla(\phi_\star A)$ has polynomial
growth too. Now, by using the Leibniz identity and an 
induction argument, one gets (i).\\ 
For (ii), we first look for a positive number $N$ such that
\begin{equation}\label{A}
\int_{\a_o}||a||^{-N}|(\phi^{-1})^{\star}u(a)|\, da<\infty\quad 
(u\in\s(\a)),
\end{equation}
where $\a_{o}$ is the complement of some compact neighborhood of the 
origin in $\a$. A change of variables following $a\leftarrow\phi(a)$ 
leads us to 

\begin{equation}\label{B}
\int_{\a'_o}||\phi(a)||^{-N}|\mbox{Jac}_{\phi}(a)|\,|u(a)|\, da
\end{equation}
where $\mbox{Jac}_{\phi}(a)$ denotes the determinant of the 
differential of $\phi$ at point $a$, and where $\a'_o$ is of the same 
type as $\a_o$ (note that the origin is fixed by the diffeomorphism 
$\phi$). With the notations adopted in the proof of \lref{GROWTH}, one 
observes that $\mbox{Jac}_{\phi}(a)$ is proportional to 
$\Pi_{\lambda\in\Phi}|\cosh(\lambda(a))|^{\dim(V_{\lambda})}$.  
Therefore, for some constant $c>0$, one has  

$$ |\mbox{Jac}_{\phi}(a)|\leq 
c\,\max\left\{ 1, 
\left(\Pi_{\lambda\in\Phi}|
\cosh(\lambda(a))|^{\dim(V_{\lambda})}\right)^{2}\right\}
= c\,\max\left\{ 1, 
\Pi_{\lambda\in\Phi}|\sinh^{2}(\lambda(a))+1|^{\dim(V_{\lambda})}\right\}.
$$  
Thus
$$
|\mbox{Jac}_{\phi}(a)|\leq 
c\,\max\left\{ 1, 
\Pi_{\lambda\in\Phi}\left(
|\sinh^{2}(\lambda(a))|+1\right)^{\dim(V_{\lambda})}\right\}
\leq c(1+||\sinh(a)||^{2})^{K}=c(1+||\phi(a)||^{2})^{K}
$$ 
for some 
$K>0$ (cf.~proof of \pref{GROWTH-P}). This last expression being lower 
than $||\phi(a)||^{N}$ for some $N>0$ and $||a||$ large enough, this provides 
the desired $N$ in order to bound (\ref{B}) hence (\ref{A}).\\
For derivatives of $(\phi^{-1})^{\star}u$, an argument as in (i) 
leads us to consider $\phi^{-1}_{\star\phi(a)}(A)\quad(a,A\in\a)$ 
that is to consider the inverse matrix $\left[ \phi_{\star 
a}\right]^{-1}$ i.e. $\cosh(a)^{-1}$. Hence, by use of \lref{COSECH}, 
one now gets

$$
\int_{\a_o}\frac{1}{||a||^{N}}|(D_{A}(\phi^{-1})^{\star}u)(a)|\,da<\infty\quad\forall 
A\in\a; u\in\s(\a).
$$
Similarly to (i), an induction yields (ii).
\EPf

Set $\pb=\a\times\a$. Identifying $\EK^\star$ with $\a$ via the symplectic 
structure $\Omega$, one can consider the partial Fourier transform~:

$$
\s'({\cal P})\stackrel{F}{\to}\s'(\pb)
$$
formally given by 

$$
Fu(a,\alpha)=\hat{u} (a,\alpha) = 
\int_{\EK}e^{{-i\Omega(\alpha,l)}}u(a,l)\, dl .
$$
One denotes its inverse by $\s'({\cal P})\stackrel{F^{-1}}{\to}\s'(\pb)$.\\
Now, led by ~\sref{REPTH}, we make the following definition.

\begin{dfn}\label{PHI}
For all $\hbar>0$, we denote by 
$\varphi_\hbar:\pb\to\pb$ the diffeomorphism defined by 
$$
\varphi_\hbar(a,\alpha)=(a,\frac{2}{\hbar}\phi(\frac{\hbar}{2}\alpha)).
$$
We denote by $\s({\cal P})\stackrel{\tau_\hbar}{\to}\s'({\cal P})$ the map 
$$
\tau_\hbar=F^{-1}\circ{\varphi_\hbar^{-1}}^{\star}\circ F.
$$ 
We set 
$$
{\cal E}_\hbar\stackrel{\mbox{def.}}{=}\tau_\hbar(\s({\cal P}))\subset\s'({\cal P}),
$$ 
and define a map ${\cal E}_{\hbar}\stackrel{T_\hbar}{\to}\s({\cal P})$ by 
$$
T_\hbar=F^{-1}\circ\varphi_\hbar^{\star}\circ F.
$$
\end{dfn}

\begin{rmk}{\rm 
In the nilpotent case  (i.e. $\g$ nilpotent), one has ${\cal E}_\hbar=\s({\cal P})$.
}\end{rmk}

\begin{prop}
\hspace{-.2cm}
\begin{enumerate}
\item[(i)] $\s({\cal P})\subset {\cal E}_\hbar$.
\item[(ii)] $T_\hbar\circ \tau_\hbar=id_{\s({\cal P})}$.
\item[(iii)] $\tau_\hbar\circ T_\hbar|_{\s({\cal P})} = id_{\s({\cal 
P})}$.
\item[(iv)] Let $\star^0_\hbar$ be the Weyl product on $\s ({\cal 
P})$ (see formula~(\ref{WEYL}) in \sref{WEYLR}). Then, the expression 
$$
a\star_\hbar b\stackrel{\mbox{def.}}{=}\tau_\hbar\left(T_\hbar a \star^0_\hbar T_\hbar 
b\right)\qquad a,b\in{\cal E}_\hbar
$$
defines an associative algebra structure on ${\cal E}_\hbar$.
\end{enumerate}
\end{prop}

\Pf By \pref{STAB}, for all $u\in \s({\cal P})$, one has 
$(F^{-1}\circ\varphi_\hbar^{\star}\circ  F)u\in \s({\cal P})$.  Hence 
$\tau_\hbar(F^{-1}\circ\varphi_\hbar^{\star}\circ F)u=u\in 
\tau_\hbar(\s({\cal P}))={\cal E}_\hbar$. Items (i), (ii) and (iii) 
follow. Associativity of the Weyl product on $\s(\cal P)$ yields (iv).
\EPf

\begin{ntt}\label{NOTATION}{\rm
Since $\Omega(\frac{2}{\hbar}\phi(\frac{\hbar}{2}\alpha),l)=
\frac{2}{\hbar}\xi(\sinh(\frac{\hbar}{2}\alpha)l)$, by setting 
$\varphi_0=id_{\pb}$ and $T_0=id_{\s(\p)}$, one gets a (separately 
continuous) map $[0,\infty)\times \s(\p)\to\s(\p)~: (\hbar,u)\to 
T_\hbar(u)$. Furthermore, a computation shows that $d_\lambda \circ F = 
\frac{1}{\lambda^n} F \circ d_\lambda$, where $d_\lambda$ is defined 
as follows~: 
$$
d_\lambda~:C^\infty(\p)\to C^\infty(\p)~: (d_\lambda 
u)(a,l)=u(a,\lambda l)\mbox{ for }\lambda\in\R_0.
$$
Therefore, one has 
$$
T_\hbar=d_{\frac{2}{\hbar}}\circ T_2\circ 
d_{\frac{\hbar}{2}}\quad(\hbar>0).
$$
In particular, ${\cal E}_\hbar=d_{\frac{2}{\hbar}}({\cal E})$, with 
${\cal E}=\tau_2(\s(\p))$. This leads to set ${\cal E}_0=\s(\p)$ and 
$u\star_0v=uv\quad\forall u,v\in\s(\p)$.
}\end{ntt}

\begin{thm}\label{WKBTHM}
Let $(M,\omega,s)$ be an elementary solvable symplectic symmetric 
space admitting a globally defined midpoint map. Let $G$ be its 
transvection group. Let $\phi~:\a\to\a$ be the twisting 
diffeomorphism defined in \pref{TWISTDIFFEO}. Then, the family 
$\{\e_\hbar \}_{\hbar\geq 0}$ (cf.~\dref{PHI}), defines a $G$-invariant 
WKB-quantization of $(M,\omega,s)$ (cf.~\dref{WKB}). More precisely, 
let $S\in C^\infty(M\times M\times M,\R)$ be defined by 
$$
S((a_1,l_1),(a_2,l_2),(a_3,l_3))=\xi\left(\oint_{1,2,3}\sinh(a_1-a_2)l_3\right)
$$
(cf.~\pref{PHASESET}). Then, for all $u$ and $v$ in $\d(M)\subset 
\e_\hbar$, the product reads~: 

\begin{equation}\label{POMME}
u\star_\hbar v(x)=\frac{1}{\hbar^{2n}}\int_{M\times M} 
e^{\frac{2i}{\hbar}S(x,x_{1},x_{2})}\, 
|\det\left(\cosh(a_{2}-a_{1})|_{\EK}\right)|\, 
u(x_{1})v(x_{2}) dx_{1} \, dx_{2}
\end{equation}
with $x_{i}=(a_{i},l_{i})\in M={\cal P}$ ($i=1,2$), and where $dx$ stands 
for the symplectic measure on $M=\p$. The phase function $S$ as well as 
the amplitude $|\det\left(\cosh(a_{2}-a_{1})|_{\EK}\right)|$ are invariant 
under the symmetries $\{s_{x}\}_{x\in M}$. 
\end{thm}

\Pf
For the sake of simplicity, we establish the product formula (\ref{POMME}) 
for $\hbar=2$. We set $\star^0=\star^0_2$ (Weyl's product) $T=T_2$,
$\tau=\tau_2$ and $\varphi=\varphi_{2}$. 
Also, in the formulae that follow, integration is taken over every 
variable $a_i,l_i,\alpha_i$ with $i=1,2$.
Let $u,v\in \d(M)$ be compactly supported. 
Then, formally, one has~:
$$
(Tu\star^{0}Tv)(a_{0},l_{0})=\int 
e^{iS^{0}((a_{0},l_{0}),(a_{1},l_{1}),(a_{2},l_{2}))}
e^{i\Omega(\alpha_{1},l_{1})} 
(\varphi^{\star}\hat{u})(a_{1},\alpha_{1})\, 
e^{i\Omega(\alpha_{2},l_{2})} 
(\varphi^{\star}\hat{v})(a_{2},\alpha_{2})
$$

$$
=\int e^{i\left[ 
\Omega(a_{2}-a_{1},l_{0})+\Omega(a_{0}-a_{2}+\alpha_{1},l_{1})
+\Omega(a_{1}-a_{0}+\alpha_{2},l_{2})
\right]}(\varphi^{\star}\hat{u})(a_{1},\alpha_{1})\, 
(\varphi^{\star}\hat{v})(a_{2},\alpha_{2})
$$
(using the definition of the Weyl product)
$$
=\int \, e^{i\Omega(a_{2}-a_{1},l_{0})}\left[ 
\int \, e^{i\Omega(a_{0}-a_{2}+\alpha_{1},l_{1})}
(\varphi^{\star}\hat{u})(a_{1},\alpha_{1}) 
\int \, e^{i\Omega(a_{1}-a_{0}+\alpha_{2},l_{2})}
(\varphi^{\star}\hat{v})(a_{2},\alpha_{2})  \right]
$$

$$
=\int \, 
e^{i\Omega(a_{2}-a_{1},l_{0})}(\varphi^{\star}\hat{u})(a_{1},a_{2}-a_{0}) 
(\varphi^{\star}\hat{v})(a_{2},a_{0}-a_{1}). 
$$
Moreover,
$$
\tau u (a,l)=\int 
e^{i\Omega(\alpha,l)}({\varphi^{-1}}^{\star}\hat{u})(a,\alpha) 
d\alpha=\int e^{i\Omega(\alpha,l)}\hat{u}(a,\phi^{-1}(\alpha))\,d\alpha=
$$

$$
=\int 
e^{i(\Omega(\alpha,l)-\Omega(\phi^{-1}(\alpha),\lambda))}u(a,\lambda)\,d\alpha 
d\lambda.
$$
Hence
$$
\tau(Tu\star^0 Tv)(a_0,l_0)=\int
e^{i(\Omega(\alpha,l_0)-\Omega(\phi^{-1}(\alpha),\lambda))}e^{i\Omega(a_2-a_1,\lambda)}
(\varphi^{\star}\hat{u})(a_{1},a_{2}-a_{0}) 
(\varphi^{\star}\hat{v})(a_{2},a_{0}-a_{1})\, d\lambda\,d\alpha=
$$

$$
\int
e^{i\left[\Omega(\alpha,l_0)-\Omega(\phi^{-1}(\alpha),\lambda)+
\Omega(a_2-a_1,\lambda) - \Omega(\phi(a_2-a_0),l_1)
-\Omega(\phi(a_0-a_1),l_2)\right]}u(a_1,l_1)\, v(a_2,l_2)\, d\lambda\,d\alpha=
$$

$$
\int
e^{i\left[\Omega(\phi(\alpha),l_0)-\Omega(\alpha,\lambda)+\Omega(a_2-a_1,\lambda)
- \Omega(\phi(a_2-a_0),l_1)
- \Omega(\phi(a_0-a_1),l_2)\right]}|Jac_{\phi}(\alpha)|\,u(a_1,l_1)\, 
v(a_2,l_2)\, d\lambda\,d\alpha
$$
(after changing the variables following
$\alpha\leftarrow\phi(\alpha)$)
$$
=\int
e^{i\left[\Omega(\phi(\alpha),l_0)-\Omega(\phi(a_2-a_0),l_1)-
\Omega(\phi(a_0-a_1),l_2)-
\Omega(\alpha-a_2+a_1,\lambda) 
\right]}|Jac_{\phi}(\alpha)|\,u(a_1,l_1)\, v(a_2,l_2)\, d\lambda\,d\alpha=
$$

$$
\int
e^{i\left[\Omega(\phi(a_2-a_1),l_0)-\Omega(\phi(a_2-a_0),l_1)-
\Omega(\phi(a_0-a_1),l_2)\right]}
|Jac_{\phi}(a_2-a_1)|\,u(a_1,l_1)\, v(a_2,l_2)=
$$

$$
\int e^{iS(x_0,x_1,x_2)}|Jac_{\phi}(a_2-a_1)|\,u(x_1)\,v(x_2)
$$
(using the definition of $\phi$). We now get the announced formula 
using item (ii) of \rref{JAC}. Dirac's condition is implied by Fedoriuk's formula 
(formula (1.5) p. 30 in \cite{FM}) up to order one in the parameter 
$\hbar$, provided $(u\star_{\hbar}v)(x)$ is interpreted as the 
oscillatory integral~: 
$$
\int\varphi(X)\exp\left(\frac{i}{\hbar}\Sigma(X)\right)\,dX
$$
with 
$$\begin{array}{lcl}
X & = & (x_{1},x_{2})\in M\times M, \\
\varphi(X) & = & |\det\left(\cosh(a_{2}-a_{1})|_{\EK}\right)|\, 
u\otimes v(X) \mbox{ and},\\
\Sigma(X) & = & 2S(x,X).
\end{array}$$
The critical point analysis of the function $\Sigma$ then tells us 
that $\Sigma$ has an isolated critical point at $X^{0}=(x,x)$. 
Moreover, the Hessian matrix $\partial^{2}_{X}\Sigma(X^{0})$ is 
proportional to 
$$
\left(
\begin{array}{cc}
0 & \Omega\\
\Omega & 0
\end{array}
\right)
$$
on $T_{X^{0}}(M\times M)=\p\times\p$. Hence 
$(u\star_{\hbar}v)(x)$ admits the desired asymptotic expansion.  
\EPf

\section{Topological algebras}\label{TA}

In this section, we analyze some topological properties of the algebras 
${\cal E}_\hbar\quad(\hbar>0)$ defined in \sref{WKBQ}. Each of these 
algebras being isomorphic, via a ``dilation", to the algebra $\e = 
\tau_2(\s(\p))$ (cf.~\nref{NOTATION}), we will, in this section, drop the 
symbol ``$\hbar$" in our discussion. As previously, we set 
$\tau=\tau_2$, $T=T_2$ and $\star=\star_2$. Moreover, we set 
$\s=\s(\p)$ and $\d=\d(M)$.

\begin{dfn}
On $\s$, let us denote by $(\, , \,)_\s$ the canonical $L^2$-inner product
$$
(u,v)_{\s}  \int u\overline{v}.
$$
Via the linear bijection $T~:\e\to\s$, one define on $\e$ 
the following inner product
$$
(a,b)_{\cal E} \stackrel{\mbox{def.}}{=} (Ta,Tb)_{\s}.
$$
The space $\e$ then becomes a pre-Hilbert space whose Hilbert completion 
is denoted by $\h$.
\end{dfn}

\begin{prop}\label{PROPTA}
\begin{enumerate}
\item[(i)] The inclusions
$$
\d\subset\s\subset\e\subset\h
$$
are dense.
\item[(ii)] For all $x\in M$ and $u,v\in\d$, one has
$$
(s_x^\star u,s_x^\star v)_{\cal E}=(u,v)_{\cal E}.
$$
In particular, the action of the transvection group $G$ on $\d$ extends
as an (unitary) action of $G$ on $\h$.
\item[(iii)] The algebra structure $\star$ on $\e$ extends to $\h$. Endowed with the 
extended product again denoted by $\star$, the space $\h$ becomes an associative topological algebra.
\item[(iv)] The group $G$ acts on $(\h,\star)$ by algebra automorphisms.
\end{enumerate}
\end{prop}

\Pf
Observe that the space $F^{-1}\d$ is stable by both transformations 
$T$ and $\tau$. Moreover, the inclusion $F^{-1}\d\subset\s$ is dense 
with respect to the $L^{2}$-topology. Hence $\tau F^{-1}\d=F^{-1}\d$ is dense 
in $\e$ and so is $\s$. 
Moreover, the map $T|_\s:\s\to\s$ is continuous with respect to the 
$L^{2}$-topology. Indeed, this follows from the fact that
\begin{eqnarray*}
|\mbox{Jac}_{\phi^{-1}}(a)|=\frac{1}{|\det(\phi_{\star_{\phi^{-1}(a)}})|}\\
=\frac{1}{|\det(\cosh(\phi^{-1}(a))|_\EK)|}\\
=\frac{1}{|\det(\cosh(\phi^{-1}(a))|_\EK)^2|^{\frac{1}{2}}}\\
=\frac{1}{|\det({\bf 1}+(\rho(a)^2)|_\EK)|^{\frac{1}{2}}}.
\end{eqnarray*}
Hence the inclusion $\s\stackrel{i}{\to}\e$, $i=\tau\circ T|_\s$, is 
continuous when $\s$ is endowed with the $L^{2}$-topology.
This yields (i).\\
One has, with $x=(a_0,l_0)$~:
{\small
$$
\begin{array}{c}
(s_x^\star u,s_x^\star v)_{\cal E} =  (F^{-1}\varphi^\star  Fs_x^\star
u,F^{-1}\varphi^\star  Fs_x^\star v)_{\s(\p)}
=(\varphi^\star  Fs_x^\star u,\varphi^\star  Fs_x^\star v)_{\s(\overline{\p})}=
\int |Jac_{\phi^{-1}}|Fs_x^\star u \overline{Fs_x^\star v}=\\
 \int \left[ |Jac_{\phi^{-1}}(\alpha)|
\int e^{i\Omega(\alpha,l)}u(2a_0-a,2 \cosh(a_0-a)l_0-l)\,dl 
\int^{-i\Omega(\alpha,l')}\overline{v} (2a_0-a, 2 \cosh(a_0 - 
a)l_0-l')\, dl'\right] da\,d\alpha =\\   
\int \left[ |Jac_{\phi^{-1}}(\alpha)|\int e^{-i
\Omega(\alpha,\lambda-2\cosh(a_0-a)l_0)} u(2a_0-a,\lambda)\, d\lambda 
\int e^{i \Omega(\alpha,\lambda'-2\cosh(a_0-a)l_0)}
\overline{v}(2a_0-a,\lambda')\, d\lambda'\right] da\,d\alpha.
\end{array}
$$
}
Since the terms $\Omega(\alpha,2\cosh(a_0-a)l_0)$ in the exponentials cancel 
each other, the latter expression equals
$$
\int da\,d\alpha |Jac_{\phi^{-1}}(\alpha)| \,\overline{F\overline{u}} 
(2a_0-a,\alpha) F \overline{v}(2a_0-a,\alpha)=(u,v)_{\cal E}.
$$
This proves (ii) while item (iii) follows from 
\cite{H}. Now, observe that for $a=\lim a_n\in\h$ with $a_n\in\d$, 
$u\in\d$ and $g\in G$, one has
\begin{eqnarray*}
ga\star gu = (\lim ga_n)\star gu \quad (\mbox{by (ii)})\\
	= \lim(ga_n\star gu) \quad (\mbox{by (iii)})\\
	= \lim g(a_n\star u) \quad(\mbox{by formula~(\ref{POMME})})\\
	= g\lim(a_n\star u) \quad (\mbox{by (ii)})\\
	= g(a\star u) \quad (\mbox{by (iii)}).
\end{eqnarray*}
This implies (iv).
\EPf
We now follow a standard procedure (\cite{H}, \cite{R1}).
For all $a\in\h$, the left multiplication $L_a:\h\to\h:b\to a\star b$ is 
a bounded operator. This yields an algebra homomorphism 
$$
\h\rightarrow\b(\h):a\to L_a
$$
into the $C^\star$-algebra $\b(\h)$ of bounded operators on $\h$.
This homomorphism is continuous. It is also injective. Indeed, 
denoting Weyl's product by $\star^0$, our Hilbert algebra $(\h,\star)$ 
is isomorphic to $(L^2(\p),\star^0)$. Denoting by $\mbox{tr}^{0}$ 
the canonical trace for Weyl's quantization (i.e 
$\mbox{tr}^{0}(u)=\int u$), one has 
$$
\mbox{tr}^{0}(a\star^0\overline{b})=(a,b)_{L^2}\quad\forall a,b\in L^2(\p).
$$
Hence, for $a\in\h$ such that $L_a=0$, one gets $(a,\h)=0$, that is $a=0$.
One therefore gets a new norm   
on $\h$~:
$$
||a||_2\stackrel{\mbox{def.}}{=}||L_a||_{\b(\h)}.
$$

\begin{prop}\label{CSTAR}
\begin{enumerate}
\item[(i)] The complex conjugation on $\d$ extends continuously to the Hilbert 
algebra $\h$ as an involution for the product $\star$. We denote this 
involution by $a\to a^*$. 
\item[(ii)] The quadruple $(\h,\star,||\,||_2,{}^*)$ is then a
(pre)-$C^*$-algebra on which the transvection group acts by 
$C^*$-algebra automorphisms.  
\end{enumerate}
\end{prop}

\Pf
Let $u\in\d$ and note that, from  $\varphi^\star F\circ F(u)=F\circ F(\varphi^\star u)$, 
one gets $T(\overline{u})=\overline{T(u)}$. 
Therefore, one obtains the pre-$C^\star$-algebra structure on $(\h,\star)$ 
by transporting the one on $(L^2(\p),\star^0)$ \cite{H} via the isomorphism 
$$
T:\h\to (L^2(\p).
$$
The rest follows from Proposition~\ref{PROPTA}.
\EPf

\section{Remarks for further developments}\label{CONCL}

A possible extension of this work is to define curved ``quantum symmetric 
solv-manifolds" analogue to flat quantum tori. That is, strict quantizations 
(in Rieffel's sense) of compact quotients of solvable symmetric spaces. 
The presence of curvature should yield interesting continuous fields of 
$C^\star$-algebras. Indeed, already at the level of the universal covering, 
which is considered in this present work, one can observe the non closeness 
of the Schwartz space $\s$ under our deformed products~: the one parameter 
``equivalence" $T_\hbar$ really ``turns" $\s$ inside $C_\infty(M)$ 
(cf.~\sref{WKBQ}). It is not even clear if there actually exists any reasonable 
Poisson subspace stable, through the deformation, under both classical 
and deformed products. In order to attend these questions, one would need 
to focuse the following points.

\begin{enumerate}
\item Establish continuous fields. That is, firstly, describe the strong 
deformation of $\s$ ($\hbar=0$) with generic fiber $\h$ (cf.~\pref{CSTAR}) 
arising from our construction. Secondly, extend this to (smooth) bounded 
functions---this is necessary in order to consider quotients later on. 
Observe that, regarding this last point, the proofs of \pref{GROWTH-P} 
and \pref{STAB} indicate that our setting should extend to bounded functions 
without difficulty just as in the flat case (cf.~\cite{R1}). 
\item Extend the present work to the whole class of solvable symplectic 
symmetric spaces which admit globally defined midpoint maps. Proposition 
3.2 in \cite{Bi1} indicates that this question should follow from an induction 
on the successive split extensions by Abelian algebras which eventually 
yields any solvable symmetric space from the data of an elementary one.
\item Study cocompact actions of discrete subgroups of automorphism groups of 
solvable symmetric spaces.
\end{enumerate}
In this context of solv-manifolds, one can then hope to investigate 
the problem of defining a quantum analogue to the Anosov property for 
classical flows.

The -approach to quantization of symmetric spaces is also aiming at 
attending quantum Riemann surfaces in a ``universal" setting, that is without 
referring to particular Hilbert space representations. This problem 
lies in the semisimple world which, from the geometric point of view as 
well as from the point of view of star representation 
theory, is  more complicated than the solvable situation considered here. 
This problem has actually been investigated by Berezin in \cite{B2}.
Since then, numerous of important works have emerged concerning this question 
(see e.g.~\cite{CGR}, \cite{Ra} and \cite{UU}) in the framework 
of Berezin-Toeplitz quantization. Comparisons between the WKB-quantization approach 
and Berezin's quantization of the hyperbolic plane have first been investigated 
by Weinstein and Qian (see \cite{W1} and \cite{Q}). Explicit computations of invariant  
admissible phase functions on the hyperbolic plane have been performed 
by Weinstein, Qian and the author. Also, some of the steps of the present 
construction pass to the case of the hyperbolic plane. For instance 
explicit Darboux charts in which the Moyal star product is $SL_2(\R)$-covariant 
have been found by the author. But, it is unclear how to define an 
intertwiner analogous to operator $\z_\hbar$ (cf.~\sref{REPTH}) from the 
data of the star representation cocycle and the twisting map 
(cf.~\dref{TWISTINGMAP}) in the case of the hyperbolic plane.

It seems also interesting to compare our framework with Fedosov's 
invariant quantization. At the formal level, Fedosov's quantization gives us a way to 
construct invariant star products on affine symplectic 
manifolds (``invariant'' means that the star product is preserved 
under the affine symplectic transformations) \cite{F1}. In the case of 
a symplectic symmetric space, such an invariant $\star$-product is essentially 
unique (see \cite{BBG}). The present work therefore deals with the problem of 
finding, in our framework of solvable symmetric spaces, ``oscillatory
integral formulae'' whose expansions are Fedosov's series.

At last, in \cite{UU}, Unterberger and Upmeier study a pseudo-differential calculus 
(called Fuchs calculus) on symmetric  
cones. On a given symmetric cone $C$, Fuchs calculus is equivariant under the
action of a solvable Lie group $G_c$ which is	the contraction of
the automorphism group of the complex tube domain $\Pi$ over $C$. 
Explicit integral formulae for the composition products of Fuchs
symbols have been obtained in \cite{U}. In some cases, our 
framework overlaps Unterberger's one and, in such cases, the formulae
we obtain for the WKB-quantizations coincide, up to some diffeomorphism 
which can be interpreted in geometric terms, with Unterberger's composition 
formulae. Despite the fact that, in our case, no symbol-operator 
correspondence is used, our work therefore appears to be closely related 
to Fuchs calculus. A deeper study of the relation with Fuchs calculus will 
be investigated in a forthcoming paper \cite{BM}.

\end{document}